\documentclass{article}

\usepackage{multicol}        
\usepackage[bottom]{footmisc}
\usepackage{amssymb}
\usepackage{amsmath}
\usepackage{mathrsfs}
\usepackage{color}
\usepackage{caption}
\usepackage[normalem]{ulem}
\usepackage{calligra}
\usepackage[T1]{fontenc}
\usepackage{lipsum}

\usepackage{epsfig}
\usepackage{graphicx}
\usepackage{placeins}
\usepackage{subfigure}

\newtheorem{theorem}{Theorem}
\newtheorem{lemma}{Lemma}

\newtheorem{definition}{Definition}

\newcommand{\be}{\begin{equation}}
\newcommand{\ee}{\end{equation}}
\newcommand{\bea}{\begin{eqnarray}}
\newcommand{\eea}{\end{eqnarray}}
\newcommand{\bean}{\begin{eqnarray*}}
\newcommand{\eean}{\end{eqnarray*}}
\newcommand{\la}{\label}

\newcommand\blfootnote[1]{%
  \begingroup
  \renewcommand\thefootnote{}\footnote{#1}%
  \addtocounter{footnote}{-1}%
  \endgroup
}

 \newcommand{\nl}{\newline}
 
\newcommand{\N}{{\mathbb{N}}}

\newcommand{\R}{{\mathbb{R}}}

 \newcommand{\cC}{{\cal C}}

 \newcommand{\cA}{{\cal A}}

\newcommand{\loc}{\rm loc}
\newcommand{\nn}{\nonumber}

\newcommand{\ino}{\int_{\xO}}

 \newcommand{\diver}{{\rm div}}

 \newcommand{\darr}[4]{{\left\{\begin{array}{ll}
   {#1}&{#2}\\[0.2cm]
   {#3}&{#4}
 \end{array}\right.}}

\newcommand{\ia}{({\rm i})}
\newcommand{\ib}{({\rm ii})}

\newcommand{\bT}{{\bf T}}

\newcommand{\finedim}{{\hfill $\Box$}}

\newcommand{\ccD}{{\mathscr D}}
\newcommand{\ccC}{{\mathscr C}}
\newcommand{\ccA}{{\mathscr A}}

\newcommand{\xa}{\alpha}

\newcommand{\xk}{\kappa}
\newcommand{\xl}{\lambda}

\newcommand{\xs}{\sigma}
\newcommand{\xS}{\Sigma}

\newcommand{\xO}{\Omega}

\title{Sharp Hardy and Hardy--Sobolev inequalities with point  singularities  on the  boundary }
\date{}

\author{
G. Barbatis\thanks{Department of Mathematics,
 National and Kapodistrian University of Athens,  15784 Athens, Greece}
 \and S. Filippas  
\thanks{Department of Mathematics and Applied Mathematics,
 University of Crete, 70013 Heraklion, Greece }
\and A. Tertikas
 \footnotemark[2]
}

\begin{document}


\maketitle
\blfootnote{Email addresses: gbarbatis@math.uoa.gr; filippas@uoc.gr; tertikas@uoc.gr}
\blfootnote{Corresponding author: A. Tertikas}

\tableofcontents

\begin{abstract}
\noindent
We  study the Hardy inequality  
when the singularity is placed on the boundary of a bounded domain in $\R^n$
that satisfies both an interior and exterior ball condition at the singularity.
We obtain the sharp Hardy constant $n^2/4$ in  case the  exterior ball is large enough and show  the necessity of the large exterior ball condition.
We improve Hardy inequality with the best constant by adding a sharp Sobolev term. We next produce criteria that lead to  characterizing  maximal potentials that improve Hardy inequality. Breaking the criteria one  produces successive improvements with sharp constants. Our approach goes through in less regular domains, like cones. In the case of a cone, contrary to the smooth case, the Sobolev constant does depend on the opening of the cone.

\end{abstract}

\vspace{11pt}


\vspace{6pt}

\section{Introduction  and main results}

For $n \geq 3$  Hardy inequality states,  that for any $u \in 
C^{\infty}_{c}(\R^n)$ there holds
\[
\int_{\R^n} |\nabla u|^2 dx 
\geq  \Big( \frac{n-2}{2}\Big)^2  \int_{\R^n} \frac{u^2}{|x|^2} dx \ ,
\]
where $\frac{(n-2)^2}{4}$ is the best constant. On the other hand Sobolev inequality reads as follows
\[
\int_{\xO} |\nabla u|^2 dx \geq S_{n} \bigg( \int_{\xO}  |u|^{\frac{2n}{n-2}}     dx \bigg)^{\frac{n-2}{n}}, \qquad u\in C^{\infty}_c(\Omega)
\]
where $S_n=
 \pi n(n-2) \left( \Gamma(\frac{n}{2})/ \Gamma(n) \right)^{\frac2n}$ is the best Sobolev constant for any domain $\Omega\subset\R^n$.

There are 
various improved versions of either  Hardy or Sobolev inequalities   in the case  of a bounded domain $\xO$ containing the origin see e.g
\cite{BV,VZ,Te,ACR,FT2002,BFT1,BFT2,FMoT,BN,BL}. We mention in particular the
following sharp Hardy--Sobolev inequality  from \cite{FT2002,AFT} that combines both inequalities
\[
\int_{\xO} |\nabla u|^2 dx 
\geq   \Big( \frac{n-2}{2}\Big)^2
 \int_{\xO} \frac{u^2}{|x|^2} dx \\
+ (n-2)^{-\frac{2(n-1)}{n}}  \, S_{n} \bigg( \int_{\xO}  X_1^{\frac{2n-2}{n-2}}|u|^{\frac{2n}{n-2}}     dx \bigg)^{\frac{n-2}{n}}
\]
for all $u \in C^{\infty}_{c}(\xO)$. Here $X_1=X_1(|x|/ D)$, with 
\[
X_1(t) = \frac{1}{1- \ln t}, \;\;\; t \in (0,1), \quad\quad D:=\sup_{x \in \xO}|x| \ .
\]

A natural question is what are the analogues of Hardy and Hardy--Sobolev inequalities in case the origin is on the boundary of $\xO$ instead of being in the interior. As we shall see, contrary to the previous case, the geometry of $\xO$ plays an important role. In the simplest case of the  half space $\R^n_{+}=\{(x',x_n):~x_n>0\}$, Hardy 
inequality with best constant reads (cf. \cite{N, Fil3})
\be\la{in10}
\int_{\R^n_{+}} |\nabla u|^2 dx \geq   \frac{n^2}{4}  \int_{\R^n_{+}} \frac{u^2}{|x|^2} dx, ~~~~~~~~~~\forall  u \in 
C^{\infty}_c(\R^n_+)\ .
\ee
In the more general case where the domain is a cone $\ccC$ with its vertex at the origin  the sharp Hardy inequality reads
(cf. \cite{N})
\be
\la{cone}
\int_{\ccC } |\nabla u|^2 dx \geq  \left( \Big( \frac{n-2}{4}\Big)^2  + \mu_1(\xS) \right) \int_{\ccC } \frac{u^2}{|x|^2} dx,  \qquad \forall  u \in C_c^{\infty}(\ccC ) \ ,
\ee
where $ \xS = \ccC \cap S^{n-1}$ and $\mu_1(\xS)$ is the first Dirichlet eigenvalue of the Dirichlet Laplace-Beltrami operator on $\xS$.
 
If on the other hand, the origin is on the boundary of a smooth near zero domain, then, related types of problems have been studied in \cite{GK, GR1, CL, GR2}. More precisely the following minimization problem has been considered for $0<s<2$ and $n \geq 4$,
 \be\la{13}
 \mu_{s}(\xO) = \inf_{u \in H^1_{0}(\xO)} \frac{\int_{\xO} |\nabla u|^2 dx }{\bigg(\int_{ \xO}  \frac{|u|^{\frac{2(n-s)}{n-2}}}{|x|^s}     dx \bigg)^{\frac{n-2}{n-s}}}  \ ,
 \ee
and it was established that the geometry of $\xO$ around zero plays an important role. In particular if the mean curvature at zero is negative then $ \mu_{s}(\xO) < \mu_{s}(\R^n_{+})$
 and there exists a minimizer for (\ref{13}). In the limit case $s=2$ the infimum  $ \mu_{2}(\xO)$ is the best Hardy  constant and under certain geometric assumptions on $\xO$ has been studied in \cite{Ccras1,Ccras2,Cfa,Csiam,FM1}.

  In \cite{F1} it was realized  that    the geometry  plays no role for the local best Hardy constant. That is,  
  for $r>0$ small enough if  we denote by $B_r$ the ball of radius $r$ centered at the origin, then for a smooth near zero domain $\xO$ one has
 \be\la{fa}
\int_{\xO \cap B_r} |\nabla u|^2 dx \geq   \frac{n^2}{4}  \int_{\xO \cap B_r} \frac{u^2}{|x|^2} dx, ~~~~~~~ u \in C^{\infty}_{c}(\xO \cap B_r) \ ,
\ee
which in particular implies the existence of a constant $\xl \geq 0$
such that
\bean
\xl  \int_{\xO} u^2 dx  + \int_{\xO} |\nabla u|^2 dx 
\geq  \frac{n^2}{4}  \int_{\xO} \frac{u^2}{|x|^2} dx , \quad u \in C^{\infty}_{c}(\xO) \ .
\eean
 
The first question we raise in this work  is to find a more quantitative result
 that connects the local inequality (\ref{fa}) to the global inequality in the half space (\ref{in10}). To state our first result we denote by $B_{\rho}(x_0)$ the ball of radius $\rho$ centered at $x_0$ and simply by  $B_{\rho}$ in case the ball is centered at the  origin; we also denote by $\cC A$ the complement of a set $A\subset\R^n$.

Throughout this work $\xO \subset \R^n$, $n \geq 2$, is a bounded domain with $0 \in \partial \xO$  satisfying an exterior ball condition
at zero, that is there exists a ball
\[B_{\rho}(- \rho e_n) \subset \cC \overline{\xO} \ .
\]
 We also denote 
\[
D:=\sup_{\Omega}|x|.
\]

\begin{theorem}
\label{b}
 There exists a positive  constant $\tau_n$ depending only on $n$ such that
if  the radius of the exterior ball satisfies $\rho\geq D/ \tau_n$ then
\be
\int_{\xO} |\nabla u|^2 dx \geq   \frac{n^2}{4}  \int_{\xO} \frac{u^2}{|x|^2} dx \ ,
\la{eq:ioa}
\ee
for all $u \in C^{\infty}_{c}(\xO)$. 
If in addition $\Omega$ satisfies an interior ball condition at $0$ then the constant $n^2/4$ is sharp. 
\end{theorem}

Thus in the case of a smooth (near zero) domain $\xO$, if the exterior ball at zero is large enough compared to the size of $\xO$ then the Hardy constant is $n^2/4$. If however the (largest) exterior ball is not large enough, at the end of Section \ref{sec3} we present an Example 
where the Hardy constant is smaller than $n^2/4$.

We next improve Hardy inequality by adding a Sobolev term:

\begin{theorem}
\label{c}
Let   $n \geq 3$. There exist positive constants $\sigma_n $ and $C_n$  that depend only on $n$ such that, if the radius of the exterior ball satisfies  $\rho\geq D/\xs_n$ the following holds true: 
\[
\int_{\xO} |\nabla u|^2 dx \geq   \frac{n^2}{4}  \int_{\xO} \frac{u^2}{|x|^2} dx \\
+  C_n \bigg( \int_{\xO}  X_1^{\frac{2n-2}{n-2}}|u|^{\frac{2n}{n-2}}     dx \bigg)^{\frac{n-2}{n}},
\]
for all $u \in C^{\infty}_{c}(\xO)$. Here $X_1=X_1(|x|/ 3D)$.
If in addition $\Omega$ satisfies an interior ball condition at $0$ then 
the exponent $(2n-2)/(n-2)$ of $X_1$ is sharp.
\end{theorem}

If  the radius of the exterior ball is small then there exists  a non negative constant $\xl$ (that depends on $\xO$) so that we have 
\[
\xl \int_{\xO} u^2 dx + \int_{\xO} |\nabla u|^2 dx 
\geq   \frac{n^2}{4}  \int_{\xO} \frac{u^2}{|x|^2} dx \\
+  C_n \Big( \int_{\xO}  X_1^{\frac{2n-2}{n-2}}|u|^{\frac{2n}{n-2}}     dx \Big)^{\frac{n-2}{n}}, ~~~~u \in C^{\infty}_{c}(\xO) \ ;
\]
for the precise statement see Theorem \ref{4c}.

Under the assumptions of Theorem \ref{c}, a simple application of Holder's inequality yields that for any $\xa>2$  there exists a positive constant $c(\xa, \xO)$ such  that
\[
\int_{\xO} |\nabla u|^2 dx \geq   \frac{n^2}{4}  \int_{\xO} \frac{u^2}{|x|^2} dx \\
+  c(\xa, \xO)  \int_{\xO} \frac{X_1^{\xa}}{|x|^2}u^{2} dx , ~~~~~u \in C^{\infty}_{c}(\xO) \ .
\]
If $c(\xa, \xO)$ is the best constant then this inequality cannot be further improved, see Theorem \ref{d}. On the other hand, as we shall see,  in the limiting case  $\xa=2$ the inequality  is also true, that is
\[
\int_{\xO} |\nabla u|^2 dx \geq   \frac{n^2}{4}  \int_{\xO} \frac{u^2}{|x|^2} dx \\
+  \frac14  \int_{\xO} \frac{X_1^{2}}{|x|^2}u^{2} dx , ~~~~~u \in C^{\infty}_{c}(\xO) \ ,
\]
and  the constant $1/4$ is sharp. In contrast with the case $\xa>2$ this inequality can be further improved.

This is a particular case of a more general situation where one has a non negative potential $V$  that for some $\xl$ non negative  and some sharp  positive constant $C$ the following inequality is true:
\be
\lambda\ino u^2dx +
\ino|\nabla u|^2dx \geq  \frac{n^2}{4} \ino\frac{u^2}{|x|^2}dx +C \ino Vu^2dx, \quad u \in  C^{\infty}_{c}(\xO) \ .
\la{in301}
\ee

In Section \ref{section:impr} we characterize maximal potentials, that is potentials $V$ such that (\ref{in301}) cannot be improved, with $C$ being the best constant for (\ref{in301}); such examples are the subcritical potentials, see Definition \ref{def:2}.
The main result of Section 5 is Theorem \ref{d}. We note that this description of maximal potentials is analogues to the description 
in \cite{Te,FT2002} for the interior point singularity case. 

In Section \ref{sub:log:impr} we consider the problem of successively improving Hardy inequality by critical potentials.
Before stating our result we first define the iterated logarithms (cf. \cite{FT2002})
\[
X_{k+1}(t) = X_k(X_1(t)), \qquad t\in (0,1],~~~~k=1,2,\ldots
\]
One can check that for $t \in (0,1)$ the series $ \sum_{i=1}^{\infty} X_1(t)X_{2}(t)\ldots X_i(t)$ converges
(see the proof of Lemma 6.3 in \cite{FT2002} or the Appendix in \cite{De}) and
that it is a strictly increasing function of $t$. We denote by $\kappa$ the unique $\kappa >1$ for which
\be
\sum_{i=1}^{\infty}X_1(1/\kappa)\ldots X_i(1/\kappa)=\frac{1}{4}.
\la{kappa}
\ee
We then have 
\begin{theorem}
\label{f}
 There exists $\sigma_n >0$ that depends only on $n$ such that if the radius of the exterior ball satisfies  $\rho\geq D /\xs_n$ the following holds true: 
\bean
 \int_{\xO} |\nabla u|^2 dx
&\geq &   \frac{n^2}{4}  \int_{\xO} \frac{u^2}{|x|^2} dx + \frac{1}{4} \sum_{i=1}^{\infty} \int_{\xO} \frac{u^2}{|x|^2}X_1^2\ldots X_i^2 dx   ,
\eean
for all $u \in C^{\infty}_{c}(\xO)$; here $X_i=X_i(|x|/ (3\kappa D))$.
If in addition $\Omega$ satisfies an interior ball condition at $0$ then the constants $1/4$ are sharp  at each step, that is
\[
\inf_{u \in C^{\infty}_{c}(\xO)} \frac{\int_{\xO} |\nabla u|^2 dx -
\frac{n^2}{4}  \int_{\xO} \frac{u^2}{|x|^2} dx }
{ \int_{\xO} \frac{u^2}{|x|^2}X_1^2 dx} =  \frac{1}{4} \ ,
\]
and for each $m=2,3,\ldots $
\[
\inf_{u \in C^{\infty}_{c}(\xO)} \frac{\int_{\xO} |\nabla u|^2 dx -
\frac{n^2}{4}  \int_{\xO} \frac{u^2}{|x|^2} dx - \frac{1}{4} \sum_{i=1}^{m-1} \int_{\xO} \frac{u^2}{|x|^2}X_1^2\ldots X_i^2 dx}
{ \int_{\xO} \frac{u^2}{|x|^2}X_1^2\ldots X_m^2 dx} =  \frac{1}{4} \ .
\]
\end{theorem}
We also have the Hardy-Sobolev analogue:
\begin{theorem}
\label{g}
Let  $n \geq 3$. There exist positive constants $\sigma_n $ and $C_n$  that depend only on $n$ such that, if the radius of the exterior ball satisfies  $\rho\geq D/\xs_n$ then   for any $m\in\N$ the following holds true:
\[
 \int_{\xO} |\nabla u|^2 dx  \geq   \frac{n^2}{4}  \int_{\xO} \frac{u^2}{|x|^2} dx +
 \frac{1}{4} \sum_{i=1}^m \int_{\xO} \frac{u^2}{|x|^2}X_1^2\ldots X_i^2 dx 
+  C_n \left( \int_{\xO}  (X_1\ldots X_{m+1})^{\frac{2n-2}{n-2}}|u|^{\frac{2n}{n-2}}     dx \right)^{\frac{n-2}{n}} ,
\]
for all $u \in C^{\infty}_{c}(\xO)$; here $X_i=X_i( |x| / (3\kappa D))$. 
If in addition $\Omega$ satisfies an interior ball condition at $0$ then the exponents $(2n-2)/(n-2)$ of $X_i$ are also sharp.
\end{theorem}

We then proceed to obtain a characterization for maximal potentials in the context of logarithmic improvements; see Theorem \ref{h}.

Analogues of these theorems hold true if the domain $\xO$ is a cone with vertex at zero and Section \ref{section:cone} is entirely
devoted to this. What is interesting in this case is that the Sobolev constant depends on the cone. As a typical result we mention here the following theorem that
refers to a bounded cone $\ccC_1 := \ccC \cap B_1$, the intersection of an infinite cone $\ccC$ with vertex at the origin with the unit ball $B_1$.
\begin{theorem}
Let $n\geq 3$. There exists a positive constant $C$ that depends only on $\Sigma$ such that
\be
\int_{\ccC_{1}}|\nabla u|^2dx \geq \Big( \frac{(n-2)^2}{4} +\mu_1(\Sigma)\Big) \int_{\ccC_{1}}\frac{u^2}{|x|^2}dx    + C \left( \int_{\ccC_{1}}  X_1^{\frac{2n-2}{n-2}}|u|^{\frac{2n}{n-2}}     dx \right)^{\frac{n-2}{n}},
\label{in99}
\ee
for all $u\in C^{\infty}_c(\ccC_1)$; here $X_1=X_1(|x|)$. The exponent $(2n-2)/(n-2)$ of $X_1$ is the best possible. Moreover the best constant $C$ for inequality (\ref{in99}) satisfies the estimate
\be
C \leq   C_n|\Sigma|^{\frac{2}{n}}
\la{in101}
\ee
for some positive constant $C_n$ that depends only on $n$.
In particular the best constant $C$ of inequality (\ref{in99})
cannot be taken to be independent of $\Sigma$.
\label{i}
\end{theorem}

Finally, a similar analysis goes through if one has potentials with multiple singularities on the boundary, see  Theorem \ref{th72} for one such result.

Our results  about point singularities on the boundary, are analogous to the case of interior point singularities 
see \cite{FT2002, Te, AFT}. We note however that whereas in the interior singularity case the geometry of $\xO$ is irrelevant, in this work  the curvature of the boundary introduces several technical difficulties even in the case of the plain Hardy inequality (\ref{eq:ioa}) as already noted in several recent  works see e.g. \cite{Ccras1, Ccras2, Cfa, Csiam, F1, FM1, GR1, GR2}.
 To overcome these difficulties we produce new improved inequalities in the flat case, see Lemmas
\ref{gera:lem:1.12im_ex}, \ref{gera:lem:123}, \ref{gera:lem_mlogs}
and then we use suitable conformal transformations thus obtaining sharp inequalities under the exterior ball assumption.

\section{Distance from  the vertex of a cone}
\la{section:cone}

In this section we consider the case of a finite cone and we obtain both homogeneous and nonhomogeneous improvements of the Hardy inequality (\ref{cone}). We pay particular attention to the special case where the cone is the half ball $B^{+}_{R}$. In this case the  estimates we obtain are stronger than in the case of a general cone and play a crucial role in our subsequent  analysis.

Let $\Sigma  \subset S^{n-1}$ be a domain in $S^{n-1}$ (that is a set that is open and connected in the relative topology)
with Lipschitz boundary. Let $\mu_k=\mu_k(\xS)$ be the $k$th Dirichlet eigenvalue of the Laplace-Beltrami operator on $\xS$
and let $\phi_k$ be a corresponding eigenfunction that is,
\[
\darr{- \Delta_{S^{n-1}} \phi_k(\omega)  =  \mu_k \phi_k(\omega),\;\; }{ \omega \in \xS \ ,} 
{\phi_k \big|_{\partial  \xS}   = 0  .}{}
\]
We may assume that $\{\phi_k\}$ is a complete orthonormal system in $L^2(\Sigma)$. We note that $\mu_1$ is a simple eigenvalue and we take $\phi_1$ to be positive.

We define
\[
\ccC = \{ x \in \R^n \smallsetminus \{0\}: ~~\frac{x}{|x|} \in \xS \} \, , \qquad
\ccC_{1} = \ccC  \cap B_1 = \{ x \in \R^n \smallsetminus \{0\}: ~~\frac{x}{|x|} \in \xS \, ,
|x|<1  \} \ .
\]
{\bf \em Proof of Theorem \ref{i}.} Let $u\in C^{\infty}_c(\ccC_1)$ be given and let
\[
u(x)=\sum_{k=1}^{\infty}u_{k}(r)\phi _{k}(\omega)
\]
be its decomposition into the spherical harmonics of $\Sigma$. We then have
\[
u_k(r)=\int_{\Sigma} u(x)\phi_k(\omega)dS(\omega).
\]
Let $\omega_{n-1}$ denote the surface measure of the unit sphere $S^{n-1}$. Throughout this proof for any radial function $G$ (which sometimes shall be written as $G(x)$ and sometimes as $G(r)$) we shall use the notation
\[
\int_{B_1}G(x)dx =\omega_{n-1}\int_0^1 G(r)r^{n-1}dr. 
\]
It then easily follows that
\bean
\int_{\ccC_1}|\nabla u|^2dx &=& \frac{1}{\omega_{n-1}} \int_{B_1} \Big( |\nabla u_1|^2 +\mu_1\frac{u_1^2}{|x|^2} \Big) dx 
 + \frac{1}{\omega_{n-1}}\sum_{k=2}^{\infty}\int_{B_1} \Big( |\nabla u_k|^2 +\mu_k\frac{u_k^2}{|x|^2}\Big)dx \\
&=& \frac{1}{\omega_{n-1}}\int_{B_1} \Big( |\nabla u_1|^2 +\mu_1\frac{u_1^2}{|x|^2} \Big) dx 
+ \int_{\ccC_1}|\nabla (u-u_1 \phi_1)|^2dx \, .
\eean
Moreover for any bounded radial function $G$ we have
\[
\int_{\ccC_1}Gu^2 dx  =	\frac{|\Sigma|}{\omega_{n-1}}\int_{B_1}Gu_1^2 dx  
+ \int_{\ccC_1}G(u-u_1 \phi_1)^2 dx \, .
\]
Therefore
\bean
&& \hspace{-1.5cm}\int_{\ccC_1}|\nabla u|^2dx - 
\Big( \frac{(n-2)^2}{4} +\mu_1\Big) \int_{\ccC_1}\frac{u^2}{|x|^2}dx \\
&=& \frac{1}{\omega_{n-1}} \bigg[ 
\int_{B_1} |\nabla u_1|^2dx   - \Big( \frac{n-2}{2} \Big)^2\int_{B_1}\frac{u_1^2}{|x|^2}dx \bigg]  \\
&& + \frac{1}{\omega_{n-1}}\sum_{k=2}^{\infty}\int_{B_1} \Big( |\nabla u_k|^2  +\mu_k\frac{u_k^2}{|x|^2}
-\Big( \frac{(n-2)^2}{4} +\mu_1\Big) \frac{u_k^2}{|x|^2}\Big)dx \\
&\geq&  \frac{1}{\omega_{n-1}} \bigg[ 
\int_{B_1} |\nabla u_1|^2dx   - \Big( \frac{n-2}{2} \Big)^2\int_{B_1}\frac{u_1^2}{|x|^2}dx \bigg]  \\
&& + \frac{1}{\omega_{n-1}} \cdot \frac{\mu_2-\mu_1}{ (\frac{n-2}{2})^2+ \mu_2}  
\sum_{k=2}^{\infty}\int_{B_1} \Big( |\nabla u_k|^2  +\mu_k\frac{u_k^2}{|x|^2}\Big) dx \\
&\geq& C_n \Big(\int_{B_1}X_1^{\frac{2(n-1)}{n-2}}|u_1|^{\frac{2n}{n-2}}dx\Big)^{\frac{n-2}{2}} + 
 \frac{\mu_2-\mu_1}{ (\frac{n-2}{2})^2+ \mu_2} \int_{\ccC_1} |\nabla (u-u_1\phi_1)|^2 dx \\
&\geq& \frac{C_n}{1+\mu_1}  \Big(\int_{\ccC_1}
X_1^{\frac{2(n-1)}{n-2}}|u_1\phi_1|^{\frac{2n}{n-2}}dx\Big)^{\frac{n-2}{2}} + 
 \frac{(\mu_2-\mu_1)S_n}{ (\frac{n-2}{2})^2+ \mu_2}
\Big(\int_{\ccC_1} |u-u_1\phi_1|^{\frac{2n}{n-2}} dx \Big)^{\frac{n-2}{2}}\\
&\geq& C(n,\Sigma)\Big(\int_{\ccC_1} X_1^{\frac{2(n-1)}{n-2}}|u|^{\frac{2n}{n-2}}dx\Big)^{\frac{n-2}{2}}.
\eean
For the optimality of the exponent, suppose to the contrary that there exists $p<(2n-2)/(n-2)$ such that
\[
\int_{\ccC_1}|\nabla u|^2dx \geq \Big( \frac{(n-2)^2}{4} +\mu_1(\Sigma)\Big) \int_{\ccC_1}\frac{u^2}{|x|^2}dx    + C \left( \int_{\ccC_1}  X_1^{p}|u|^{\frac{2n}{n-2}}     dx \right)^{\frac{n-2}{n}},
\]
for all $u\in C^{\infty}_c(\ccC_1)$. Considering functions $u$ of the form $u(x)=v(r)\phi_1(\omega)$ with
$v(1)=0$ we obtain that any such $v$ satisfies
\be
\int_{B_1}|\nabla v|^2 dx \geq \Big( \frac{n-2}{2}\Big)^2  \int_{B_1}\frac{v^2}{|x|^2}dx 
+    C \left( \int_{B_1}  X_1^{p}|v|^{\frac{2n}{n-2}}     dx \right)^{\frac{n-2}{n}}.
\la{123}
\ee
This is a contradiction since the best exponent of $X_1$ in (\ref{123}) is $2(n-1)/(n-2)$; see \cite{FT2002}.

To prove estimate (\ref{in101})  we test inequality (\ref{in99}) with a function of the form
$u(x)=v(r)\phi_1(\omega)$. Then an easy calculation gives
\bean
C&\leq& \frac{\int_{\ccC_1}|\nabla u|^2dx - \Big( \frac{(n-2)^2}{4} +\mu_1\Big) \int_{\ccC_1}\frac{u^2}{|x|^2}dx}
{\Big(\int_{\ccC_1}|u|^{\frac{2n}{n-2}} X_1^{\frac{2n-2}{n-2}}dx \Big)^{\frac{n-2}{n}}} \\
&=&   \frac{\frac{1}{\omega_{n-1}} \int_{B_1} \Big[ |\nabla v|^2dx - \Big( \frac{n-2}{2} \Big)^2 \frac{v^2}{|x|^2}\Big] dx}
{\Big(\frac{1}{\omega_{n-1}}\int_{B_1}|v|^{\frac{2n}{n-2}} X_1^{\frac{2n-2}{n-2}}dx \Big)^{\frac{n-2}{n}}
\Big( \int_{\Sigma}|\phi_1|^{\frac{2n}{n-2}} dS \Big)^{\frac{n-2}{n}}}
\eean
Minimizing with respect to $v$ (see \cite[Theorem B]{AFT}) we conclude that
\[
C\leq \frac{\omega_{n-1}^{-\frac{2}{n}}(n-2)^{-\frac{2}{n}}S_n }{\Big(\int_{\Sigma}|\phi_1|^{\frac{2n}{n-2}}dS\Big)^{\frac{n-2}{n}}}.
\]
By the normalization of $\phi_1$ and H\"{o}lder inequality we conclude that
\[
C\leq \omega_{n-1}^{-\frac{2(n-1)}{n}}(n-2)^{-\frac{2}{n}}S_n |\Sigma|^{\frac{2}{n}},
\]
which concludes the proof. $\hfill\Box$

In a similar fashion we obtain
\begin{theorem}
Let $n\geq 3$. There exists a constant $C$ that depends only on $\Sigma$ such that for any $m\in\N$
\bea
\int_{\ccC_{1}}|\nabla u|^2dx &\geq&
\Big( \frac{(n-2)^2}{4} +\mu_1(\Sigma)\Big)\int_{\ccC_{1}}\frac{u^2}{|x|^2}dx   + \frac{1}{4}\sum_{i=1}^m
\int_{\ccC_{1}}\frac{u^2}{|x|^2}X_1^2\ldots X_i^2 dx  \nonumber \\
&&+ C\left( \int_{\ccC_{1}}  (X_1\ldots X_{m+1})^{\frac{2n-2}{n-2}}|u|^{\frac{2n}{n-2}}     dx \right)^{\frac{n-2}{n}},
\la{gera:eq:xxxx_cone}
\eea
for all $u\in C^{\infty}_c(\ccC_1)$; here $X_i=X_i(|x|)$. Each constant $1/4$ is the best possible, that is,
\[
\inf_{u \in C^{\infty}_{c}(\ccC_1) } \frac{\int_{\ccC_1} |\nabla u|^2 dx - \big( (\frac{n-2}{2})^2 +\mu_1(\Sigma) \big) \int_{\ccC_1} \frac{u^2}{|x|^2} dx }
{ \int_{\ccC_1} \frac{u^2}{|x|^2}X_1^2 dx} =  \frac{1}{4} \ ,
\]
and for each $m=2,3,\ldots $
\[
\inf_{u \in C^{\infty}_{c}(\ccC_1)} \frac{\int_{\ccC_1} |\nabla u|^2 dx -
\big( (\frac{n-2}{2})^2 +\mu_1(\Sigma) \big)  \int_{\ccC_1} \frac{u^2}{|x|^2} dx - \frac{1}{4} \sum_{i=1}^{m-1} \int_{\ccC_1} \frac{u^2}{|x|^2}X_1^2\ldots X_i^2 dx}
{ \int_{\ccC_1} \frac{u^2}{|x|^2}X_1^2\ldots X_m^2 dx} =  \frac{1}{4}  \, .
\]
The exponent $(2n-2)/(n-2)$
is the best possible. Moreover the best constant $C$ for inequality (\ref{gera:eq:xxxx_cone}) satisfies the estimate
\be
C \leq   C_n|\Sigma|^{\frac{2}{n}}
\la{indep_log}
\ee
for some positive constant $C_n$ that depends only on $n$.
In particular the best constant $C$ of inequality (\ref{gera:eq:xxxx_cone}) cannot be taken to be independent of $\Sigma$.
\label{thm_hardy_sob_BR_cone}
\end{theorem}
{\bf \em Proof.} Arguing as in the proof of Theorem~\ref{i} we arrive at

\bean
&& \int_{\ccC_{1}}|\nabla u|^2dx - 
\Big( \frac{(n-2)^2}{4} +\mu_1\Big) \int_{\ccC_{1}}\frac{u^2}{|x|^2}dx
-\frac{1}{4}\sum_{i=1}^m\frac{u^2}{|x|^2}X_1^2\ldots X_i^2 dx \\
&=& \frac{1}{\omega_{n-1}} \bigg[ 
\int_{B_1} |\nabla u_1|^2dx   - \Big( \frac{n-2}{2} \Big)^2\int_{B_1}\frac{u_1^2}{|x|^2}dx
-\frac{1}{4}\sum_{i=1}^m\int_{B_1}\frac{u_1^2}{|x|^2}X_1^2\ldots X_i^2 dx
 \bigg]  \\
&& + \frac{1}{\omega_{n-1}}\sum_{k=2}^{\infty}\int_{B_1} \Big( |\nabla u_k|^2  +\mu_k\frac{u_k^2}{|x|^2}
-\Big( \frac{(n-2)^2}{4} +\mu_1\Big) \frac{u_k^2}{|x|^2} -\frac{1}{4}\sum_{i=1}^m 
\frac{u_k^2}{|x|^2} X_1^2\ldots X_i^2 \Big)dx \\
&\geq&  \frac{1}{\omega_{n-1}} \bigg[ 
\int_{B_1} |\nabla u_1|^2dx   - \Big( \frac{n-2}{2} \Big)^2\int_{B_1}\frac{u_1^2}{|x|^2}dx
- \frac{1}{4}\sum_{i=1}^m\int_{B_1}\frac{u_1^2}{|x|^2}X_1^2\ldots X_i^2 dx \bigg]  \\
&& + \frac{\mu_2-\mu_1}{2\big[(\frac{n-2}{2})^2+ \mu_2\big]} \frac{1}{\omega_{n-1}} 
\sum_{k=2}^{\infty}\int_{B_1} \Big( |\nabla u_k|^2  +\mu_k\frac{u_k^2}{|x|^2}\Big) dx \\
&\geq& C_n  \Big(\int_{B_1}(X_1\ldots X_{m+1})^{\frac{2(n-1)}{n-2}}|u_1|^{\frac{2n}{n-2}}dx\Big)^{\frac{n-2}{2}} + 
 \frac{\mu_2-\mu_1}{2\big[ (\frac{n-2}{2})^2+ \mu_2\big]} \int_{\ccC_1} |\nabla (u-u_1\phi_1)|^2 dx \\
&\geq& \frac{C_n}{1+\mu_1}  \Big(\int_{\ccC_1}
(X_1\ldots X_{m+1})^{\frac{2(n-1)}{n-2}}|u_1\phi_1|^{\frac{2n}{n-2}}dx\Big)^{\frac{n-2}{2}} + 
 \frac{(\mu_2-\mu_1)S_n}{2\big[ (\frac{n-2}{2})^2+ \mu_2\big]}
\Big(\int_{\ccC_1} |u-u_1\phi_1|^{\frac{2n}{n-2}} dx \Big)^{\frac{n-2}{2}}\\
&\geq& C(n,\Sigma)  \Big(\int_{\ccC_1} (X_1\ldots X_{m+1})^{\frac{2(n-1)}{n-2}}|u|^{\frac{2n}{n-2}}dx\Big)^{\frac{n-2}{2}}.
\eean
For the optimality of the constants $1/4$ we make once again the choice $u(x)=v(r)\phi_1(\omega)$ to conclude that
\bean
&& \inf_{ C^{\infty}_c(\ccC_1)}
\frac{  \int_{\ccC_1}|\nabla u|^2dx - \big( (\frac{n-2}{2})^2 +\mu_1(\Sigma) \big)
\int_{\ccC_1}\frac{u^2}{|x|^2}dx   - \frac{1}{4}\sum_{i=1}^{m-1}
\int_{\ccC_1}\frac{u^2}{|x|^2}X_1^2\ldots X_i^2 dx }{ \int_{\ccC_1}\frac{u^2}{|x|^2}X_1^2\ldots X_{m}^2
dx} \\
&\leq&  \inf_{ C^{\infty}_c(B_1)}
\frac{  \int_{B_1}|\nabla v|^2dx - \big( \frac{n-2}{2} \big)^2
\int_{B_1}\frac{v^2}{|x|^2}dx   - \frac{1}{4}\sum_{i=1}^{m-1}
\int_{B_1}\frac{v^2}{|x|^2}X_1^2\ldots X_i^2 dx }{ \int_{B_1}\frac{v^2}{|x|^2}X_1^2\ldots X_{m}^2
dx}\\
&=&\frac{1}{4},
\eean
by \cite[Theorem 6.1]{FT2002}. The optimality of the exponent in the Sobolev term follows as before from the optimality of the corresponding exponent of the Hardy-Sobolev inequality for an interior point, \cite[Theorem A']{FT2002}. 

Finally to prove estimate (\ref{indep_log}) we once again test inequality (\ref{gera:eq:xxxx_cone}) with a function of the form
$u(x)=v(r)\phi_1(\omega)$. We then obtain
\bean
C&\leq& \frac{\int_{\ccC_1}|\nabla u|^2dx - \big( \frac{(n-2)^2}{4} +\mu_1\big) \int_{\ccC_1}\frac{u^2}{|x|^2}dx
- \frac{1}{4}\sum_{i=1}^{m}
\int_{\ccC_1}\frac{u^2}{|x|^2}X_1^2\ldots X_i^2 dx }{\Big(\int_{\ccC_1}|u|^{\frac{2n}{n-2}} (X_1\ldots X_{m+1})^{\frac{2n-2}{n-2}}dx \Big)^{\frac{n-2}{n}}} \\
&=&   \frac{\frac{1}{\omega_{n-1}} \int_{B_1} \Big[ |\nabla v|^2dx - \big( \frac{n-2}{2} \big)^2 \frac{v^2}{|x|^2}
-\frac{1}{4}\sum_{i=1}^{m}\int_{B_1}\frac{v^2}{|x|^2}X_1^2\ldots X_i^2 \Big] dx}
{\Big(\frac{1}{\omega_{n-1}}\int_{B_1}|v|^{\frac{2n}{n-2}} (X_1\ldots X_{m+1})^{\frac{2n-2}{n-2}}dx \Big)^{\frac{n-2}{n}}
\Big( \int_{\Sigma}|\phi_1|^{\frac{2n}{n-2}} dS \Big)^{\frac{n-2}{n}}}
\eean
Minimizing with respect to $v$ (see \cite[Theorem B]{AFT}) and using H\"{o}lder inequality we conclude that
\[
C\leq \frac{\omega_{n-1}^{-\frac{2}{n}}(n-2)^{-\frac{2(n-1)}{n}}S_n }{\Big(\int_{\Sigma}|\phi_1|^{\frac{2n}{n-2}}dS\Big)^{\frac{n-2}{n}}}
\leq \omega_{n-1}^{-\frac{2}{n}}(n-2)^{-\frac{2}{n}}S_n |\Sigma|^{\frac{2}{n}},
\]
which concludes the proof. $\hfill\Box$

\begin{theorem}
Let $n\geq 2$. There holds
\[
\int_{\ccC_{1}}|\nabla u|^2dx \geq \bigg(   \Big( \frac{n-2}{2}\Big)^2 +\mu_1(\Sigma)\bigg)
\int_{\ccC_{1}}\frac{u^2}{|x|^2}dx   +\frac{1}{4}\sum_{i=1}^{\infty}\int_{\ccC_{1}}\frac{u^2}{|x|^2}X_1^2\ldots X_i^2dx 
\]
for all $u\in C^{\infty}_c(\ccC_{1})$; here $X_i=X_i(|x|)$. Each constant $1/4$ is sharp.
\end{theorem}
{\bf \em Proof.} This follows from Theorem \ref{thm_hardy_sob_BR_cone} by letting $m\to +\infty$. The optimality of the constants $1/4$ has been established in Theorem~\ref{thm_hardy_sob_BR_cone}.
$\hfill\Box$

\subsection{Improved Hardy inequalities in half balls}

The case of half ball where  $\Sigma=S^{n-1}_+$ is of particular importance for our approach. In this case the Hardy constant becomes
\[
 \Big( \frac{n-2}{2} \Big)^2 +n-1 =\frac{n^2}{4},
\]
and the Sobolev constants of Theorems \ref{i} and \ref{thm_hardy_sob_BR_cone}
depend only on $n$. As a special case of the previous results
we have the following sharp inequalities for all functions $u\in C^{\infty}_c(B_R^+)$:
\be
\int_{B_{R}^+}|\nabla u|^2dx \geq\frac{n^2}{4}\int_{B_{R}^+}\frac{u^2}{|x|^2}dx    + C_n \left( \int_{B_R^+}  X_1^{\frac{2n-2}{n-2}}|u|^{\frac{2n}{n-2}}     dx \right)^{\frac{n-2}{n}},
\label{gera:h_sob_br}
\ee
\bea
\int_{B_{R}^+}|\nabla u|^2dx &\geq&\frac{n^2}{4}\int_{B_{R}^+}\frac{u^2}{|x|^2}dx   + \frac{1}{4}\sum_{i=1}^m
\int_{B_R^+}\frac{u^2}{|x|^2}X_1^2\ldots X_i^2 dx  \nonumber \\
&&+ C_n \left( \int_{B_R^{+}}  (X_1\ldots X_{m+1})^{\frac{2n-2}{n-2}}|u|^{\frac{2n}{n-2}}     dx \right)^{\frac{n-2}{n}},
\la{gera:eq:xxxx}
\eea
\be
\int_{B_{R}^+}|\nabla u|^2dx \geq\frac{n^2}{4}\int_{B_{R}^+}\frac{u^2}{|x|^2}dx   +\frac{1}{4}\sum_{i=1}^{\infty}\int_{B_{R}^+}\frac{u^2}{|x|^2}X_1^2\ldots X_i^2dx  \ ,
\la{gera:eq:1.12im}
\ee
where $X_i=X_i(|x|/ R)$.

In these inequalities the singularity lies on a flat part of the boundary. However if the boundary is not flat near the singularity, then curvature plays a role.  To overcome these difficulties, in the next three lemmas we establish stronger versions of (\ref{gera:h_sob_br}), (\ref{gera:eq:xxxx}) and (\ref{gera:eq:1.12im}) that will be used to prove Theorems \ref{c}, \ref{f} and \ref{g}.

We recall (cf.(\ref{kappa})) that $\kappa$ is the unique $\kappa >1$ for which $\sum_{i=1}^{\infty}X_1(1/\kappa)\ldots X_i(1/\kappa)=\frac{1}{4}$.
We also denote for $t \in (0,1)$,
\[
\eta(t):=  \sum_{i=1}^{\infty} X_1(t)X_2(t)\ldots  X_i(t), \quad 
B(t):=  \sum_{i=1}^{\infty} X_1^2(t)X_2^2(t)\ldots  X_i^2(t) \ .
\]
Using the identity
\[
 \frac{d}{dt}X_k (t) =\frac{1}{t}X_1(t)\ldots X_{k-1}(t)X_k^2(t)
\]
we easily obtain cf \cite{BFT2}
\[
\frac{d}{dt}\eta(t) =\frac{1}{2t} \big(\eta(t)^2 +B(t) \big)  \ , \quad t \in (0,1) \ .
\]
We next have.
\begin{lemma}
For any $R>0$ there holds
\be
\int_{B_{R}^+}|\nabla u|^2dx \geq
\frac{n^2}{4}\int_{B_{R}^+}\frac{u^2}{|x|^2}dx   +\frac{1}{4}\sum_{i=1}^{\infty}\int_{B_{R}^+}\frac{u^2}{|x|^2}X_1^2\ldots X_i^2dx  + \frac{1}{8R^{1/2}}\int_{B_{R}^+}\frac{u^2}{|x|^{3/2}}dx  ,
\la{gera:eq:1.12im_ex}
\ee
for all $u\in C^{\infty}_c(B_{R}^+)$; here $X_i=X_i(|x|/\kappa R)$.
\label{gera:lem:1.12im_ex}
\end{lemma}
\noindent{\bf \em Proof.} Let ${\bf T }$ be a $C^1$ vector field in $B_R^+$ and $u\in C^{\infty}_c(B_R^+)$. We have
\[
\int_{B_R^+} \diver{\bf T}\,  u^2 dx =- 2\int_{B_R^+} u\nabla u\cdot {\bf T}dx \leq \int_{B_R^+} \big( |\nabla u|^2 + |{\bf T}|^2 u^2 \big) dx 
\]
and therefore
\[
\int_{B_R^+}|\nabla u|^2dx \geq \int_{B_R^+} \Big( \diver {\bf T} -|{\bf T}|^2 \Big)u^2 dx \; .
\]
We shall apply this for the vector field
\[
{\bf T} =\frac{n}{2}\frac{x}{|x|^2}  -\frac{e_n}{x_n} +\frac{\eta}{2}\frac{x}{|x|^2} +\frac{1}{2(R^{1/2}- |x|^{1/2})}\frac{x}{|x|^{3/2}} \, , \qquad (\eta=\eta(|x|/\kappa R)).
\]
We  have
\bean
\diver\bT&  =& \frac{n(n-2)}{2|x|^2} +\frac{1}{x_n^2} +\frac{n-2}{2|x|^2}\eta  +\frac{\eta^2 +B}{4|x|^2} \\
&& + \frac{n-\frac{3}{2}}{2(R^{1/2} -|x|^{1/2})|x|^{3/2}} + \frac{1}{4(R^{1/2} -|x|^{1/2})^2|x|} ,
\eean
hence
\bean
\diver\bT -|\bT|^2 &=&\frac{n^2}{4|x|^2} +\frac{1}{4|x|^2}\sum_{i=1}^{\infty}X_1^2\ldots X_i^2  
+\frac{\frac{1}{2}-\eta}{2(R^{1/2} -|x|^{1/2})|x|^{3/2}} \\
&\geq& \frac{n^2}{4|x|^2} +\frac{1}{4|x|^2}\sum_{i=1}^{\infty}X_1^2\ldots X_i^2  
+\frac{1}{8R^{1/2}|x|^{3/2}},
\eean
where in the last inequality we used that $\eta \leq \frac14$, because of the choice of $\xk$,
and the result follows. $\hfill\Box$

\begin{lemma}
Let $n\geq 3$. There exists a constant $C_n$ that depends only on $n$ such that for any $R>0$ there holds
\[
\int_{B_{R}^+}|\nabla u|^2dx \geq\frac{n^2}{4}\int_{B_{R}^+}\frac{u^2}{|x|^2}dx    + 
\frac{1}{16R^{1/2}}\int_{B_R^+}\frac{u^2}{|x|^{3/2}}dx 
 + C_n \left( \int_{B_R^+}  X_1^{\frac{2n-2}{n-2}}|u|^{\frac{2n}{n-2}} dx \right)^{\frac{n-2}{n}},
\]
for all $u\in C^{\infty}_c(B_R^+)$; here $X_1=X_1(|x|/R)$.
\label{gera:lem:123}
\end{lemma}
{\bf \em Proof.} The result follows by taking a convex combination of (\ref{gera:h_sob_br}) and (\ref{gera:eq:1.12im_ex}) and discarding the logarithmic terms that do not come with the sharp constant; see also the next lemma.
$\hfill\Box$

\begin{lemma}
Let $n\geq 3$ and $m\in\N$. There exists a constant $C_n$ that depends only on $n$ such that for all $R>0$ there holds
\bean
\int_{B_{R}^+}|\nabla u|^2dx &\geq&\frac{n^2}{4}\int_{B_{R}^+}\frac{u^2}{|x|^2}dx   + \frac{1}{4}\sum_{i=1}^m
\int_{B_R^+}\frac{u^2}{|x|^2}X_1^2\ldots X_i^2 dx  \\
&& \hspace{-2cm} + \frac{1}{16R^{1/2}}\int_{B_{R}^+}\frac{u^2}{|x|^{3/2}}dx 
+ C_n \left( \int_{B_R^+}  (X_1\ldots X_{m+1})^{\frac{2n-2}{n-2}}|u|^{\frac{2n}{n-2}}     dx \right)^{\frac{n-2}{n}},
\eean
for all $u\in C^{\infty}_c(B_R^+)$; here $X_i=X_i(|x|/\kappa R)$.
\label{gera:lem_mlogs}
\end{lemma}
{\bf \em Proof.} This follows by taking a convex combination of (\ref{gera:eq:xxxx}) and (\ref{gera:eq:1.12im_ex}).
$\hfill\Box$

\section{Hardy inequality in  bounded  domains}
\la{sec3}

In this section  we provide the proof of Theorem \ref{b} and use an example to establish the necessity of a relatively large exterior ball assumption.
We also analyse the Hardy constant in the case of annuli (see Theorem \ref{a}).

We initially establish that $n^2/4$ is an upper bound for the Hardy constant under an interior ball condition.
\begin{lemma}
If $\xO$  satisfies an interior ball condition at 0 then
for any $r>0$ we have
\[
\inf_{u\in C^{\infty}_c(\xO \cap B_r)} 
\frac{\int_{\xO \cap B_r} |\nabla u|^2 dx}{\int_{\xO \cap B_r} \frac{u^2}{|x|^2} dx } \leq \frac{n^2}{4}.
\]
\label{lem:sharpness}
\end{lemma}
\noindent{\bf \em Proof.} Without loss of generality we may assume that the interior ball is $B_{\rho}(\rho e_n)$ and satisfies
$B_{\rho}(\rho e_n)  \subset \Omega \cap B_r$,
therefore it is enough to establish that
\[
\inf_{u\in C^{\infty}_c(B_{\rho}(\rho e_n))} 
\frac{\int_{B_{\rho}(\rho e_n)} |\nabla u|^2 dx}{\int_{B_{\rho}(\rho e_n)} \frac{u^2}{|x|^2} dx } \leq \frac{n^2}{4}.
\]
Using a scaling argument we find that this infimum is equal to
\[
\inf_{u\in C^{\infty}_c(\R^n_+)} 
\frac{\int_{\R^n_+} |\nabla u|^2 dx}{\int_{\R^n_+} \frac{u^2}{|x|^2} dx } ,
\]
which is equal to $n^2/4$. $\hfill\Box$

We shall next prove a result about annuli. We use the notation 
 \[
\ccD( x_0 ; r_1,r_2) :=\{ x\in\R^n \; : \; r_1< |x-x_0|<r_2\}.
\]
or simply $\ccD( r_1,r_2)$ in case $x_0=0$. Also, $e_n$ shall denote the unit vector in the $x_n$ direction.

\begin{theorem}
\la{a}
Let $n\geq 2$ and let $\lambda_{\tau}$ denote the best constant for the Hardy inequality
\be
\la{in22}
\hspace{-1.1cm}  \int_{\ccD( \rho ,\rho(1+\tau))} |\nabla u|^2 dx     \geq  
 \lambda_{\tau}  \int_{\ccD( \rho ,\rho(1+\tau))} \frac{u^2}{ |x- \rho e_n|^2} dx  \, , \qquad u\in C^{\infty}_c(\ccD( \rho ,\rho(1+\tau))).
\ee
There exists a constant $\tau_n>0$ which depends only on $n$ such that 

$\quad\ia$ For all $0<\tau\leq \tau_n$ there holds $\lambda_{\tau}=n^2/4$

$\quad\ib$ For all $\tau >\tau_n$ there holds $\lambda_{\tau}<n^2/4$.

Moreover $\lambda_{\tau}$ is strictly decreasing in $(\tau_n,+\infty)$ and $\lim_{\tau\to+\infty}\lambda_{\tau}=
(n-2)^2/4$.
\end{theorem}

\noindent {\bf\em Proof.} It is enough to establish the result for $\rho=1$, the general case then follows by scaling.
To prove $\ia$ it is enough to establish that for small enough $\tau>0$ we have inequality (\ref{in22}).
We apply (\ref{gera:eq:1.12im}) with $R=2$ where we place the singularity at $e_n$ and we obtain the inequality
\be\la{hi:59}
\int_{B_1} |\nabla u|^2 dx \geq   \frac{n^2}{4}  \int_{B_1 } \frac{u^2}{|x-e_n|^2} dx  + \frac{1}{4}  \int_{B_1 } \frac{u^2}{|x-e_n|^2}X_1^2 dx,  
\quad\forall  u \in C^{\infty}_{0}(B_1)  \ , 
\ee
where $X_1=X_1(|x-e_n|/2)$. Next we apply the Kelvin transform
\[
u(x)= |y|^{n-2} v(y), \qquad y = \frac{x}{|x|^2}   \ .
\]
Then by standard  calculations using the conformality of the Kelvin transform we have
\[
\int_{B_1} |\nabla u(x)|^2 dx = \int_{\cC B_1} |\nabla v(y)|^2 dy \ ,
\]
and since 
\[
|x -e_n| = \frac{|y -e_n|}{|y|}  \ ,
\]
inequality (\ref{hi:59}) takes the equivalent form
\be
\la{hi:63}
\int_{\cC B_1} |\nabla v|^2 dy  \geq     \frac{n^2}{4}  \int_{\cC B_1 } \frac{v^2}{|y|^2 |y-e_n|^2} dy 
   + \frac{1}{4}  \int_{\cC B_1 } \frac{X_1^2}{|y|^2 |y-e_n|^2}v^2 dy  
\ee
for all $v \in C^{\infty}_c( \cC B_1)$, where $X_1=X_1( |y-e_n| / 2|y| )$.

It follows from (\ref{hi:63}) that for any $\tau>0$ and any $v \in C^{\infty}_c(B_{1+\tau} \smallsetminus B_{1})$ there holds
\[
 \int_{B_{1+\tau} \smallsetminus B_{1}} |\nabla v|^2 dy    \geq 
  \frac{n^2}{4}  \int_{B_{1+\tau} \smallsetminus B_{1}} \frac{v^2}{ |y-  e_n|^2} dy   
   +   \frac{1}{4} \int_{B_{1+\tau} \smallsetminus B_{1}} \frac{  n^2  -n^2 |y|^2 + X_1^2\left( \frac{|y- e_n|}{2|y|}  \right) }{|y|^2 |y-  e_n|^2}  v^2   dy  \ .
\]
To conclude the proof  it suffices  to show that the last term above is nonnegative for  small enough $\tau>0$. For this it is enough to have the inequality
\[
 X_1^2\Big( \frac{|y- e_n|}{2|y|}\Big) \geq n^2 \left(|y|^2-1 \right) \ ,  \qquad 1 < |y| < 1+ \tau \ .
\]
Writing $|y|=1+t$,  $0 < t < \tau$, we have that  $|y-e_n| \geq t$. Hence
\[
 X_1^2 \Big( \frac{|y- e_n|}{2|y|}\Big) \geq  X_1^2 \Big( \frac{t}{2(t+1)}\Big) \ ,
\]
and therefore it is enough to have
\[
X_1^2 \Big( \frac{t}{2(t+1)}\Big) \geq n^2 t (t+2),  \qquad 0<t<\tau \ .
\]
Since 
$\lim_{t \to 0^{+}} X_1^2(t)/t = +\infty$, the result follows.

We shall next establish that the set of all $\tau >0$ for which inequality (\ref{in22}) holds true 
is bounded and therefore we may define
\[
\tau_n= \sup\{ \tau >0 : \mbox{ inequality (\ref{in22}) holds true} \}.
\]
For this we first note that for $\tau>2$ we have the inclusion
\[
B_{\tau}\setminus B_2 \subset B_{1+\tau}(-e_n) \setminus B_1(-e_n).
\]
and therefore
\[
\inf_{C^{\infty}_c( \ccD( -e_n ; 1, 1+\tau ))} 
\frac{ \int_{\ccD( -e_n ; 1, 1+\tau )} |\nabla u|^2dx}{\int_{\ccD( -e_n ; 1, 1+\tau )}
\frac{u^2}{|x|^2}dx}
\leq
\inf_{C^{\infty}_c(B_{\tau}\setminus B_2 ) }
\frac{ \int_{B_{\tau}\setminus B_2 } |\nabla u|^2dx}{\int_{B_{\tau}\setminus B_2 } \frac{u^2}{|x|^2}dx}.
\]
Using the radial function
\[
 u(r) =r^{-\frac{n-2}{2}} \sin \Big( \frac{\ln (r/2) \, \pi}{ \ln(\tau/2)} \Big) , \quad\quad  2<r<\tau ,
\]
we easily see that the last infimum is equal to
$
\big( \frac{n-2}{2}  \big)^2 + \big( \frac{\pi}{\ln\frac{\tau}{2}}  \big)^2
$
and in particular it is smaller than $n^2/4$ if
\[
\tau >2e^{\frac{\pi}{\sqrt{n-1}}}.
\]
This implies  the existence of an $H^1_0$ minimizer (see e.g. \cite{GR1}, Theorem 4.2)  and therefore the strict monotonicity of $\lambda_{\tau}$ for $\tau>\tau_n$. The above computation also gives that
$\lim_{\tau\to +\infty}\lambda_{\tau}\leq (n-2)^2/4$; this combined with the standard Hardy inequality gives
$\lim_{\tau\to +\infty}\lambda_{\tau}= (n-2)^2/4$ thus concluding the proof of the theorem. $\hfill\Box$

We next have

{\bf\em Proof of Theorem \ref{b}:} As we shall see, the constant $\tau_n$ of Theorem \ref{b} is the same as that of Theorem \ref{a} above.
Since $\Omega \cap B_{\rho\tau_n} \subset \ccD(-\rho e_n ; \rho , \rho(1+\tau_n))$, it follows from Theorem \ref{a} that
\be
\int_{\xO \cap B_{\rho\tau_n}}  |\nabla u|^2 dx \geq   \frac{n^2}{4}  \int_{\xO \cap B_{\rho\tau_n}}  \frac{u^2}{|x|^2} dx  , \quad\quad
u\in C^{\infty}_c(\xO \cap B_{\rho\tau_n}).
\la{555b}
\ee
The assumption $\rho\tau_n \geq D$ implies $\Omega\subset B_{\rho\tau_n}$ and therefore (\ref{eq:ioa}) follows from (\ref{555b}).
The sharpness of the constant $n^2/4$ follows directly from Lemma \ref{lem:sharpness}. $\hfill\Box$

It is natural to ask whether the assumption of having a large exterior ball at zero is necessary in order to have the Hardy inequality with constant $n^2/4$. In the following example we will see that for small exterior balls inequality (\ref{eq:ioa}) fails.

%

{\bf Example.} Given $\rho \in (0,1/2)$ and $\theta\in (0,\pi/2)$ we define the domain
\[
\ccA_{\rho,\theta} =\{ x=(x',x_n)\in B_1 \; : \; x_n<\cot\theta |x'| \mbox{ and } |x-\rho e_n| >\rho \}.
\]
Let $\Omega$ be a domain containing $\ccA_{\rho,\theta}$ and having the same largest exterior ball at zero, namely $B(\rho e_n,\rho)$.

We denote by $\lambda_1(n,\theta)$ the first Dirichlet eigenvalue of the Laplace operator on the spherical cap
\[
\Sigma_{\theta}= \{ (x',x_n) \in S^{n-1} : x_n< \cot\theta \; |x'| \}.
\]
By monotonicity it follows that for $\theta <\pi/2$ we have $\lambda_1(n,\theta) < \lambda_1(n,\pi/2) =n-1$.
We shall prove that if
\be
\rho < \frac{1}{2\cos \theta} e^{- \frac{\pi}{\sqrt{n-1-\lambda_1(n,\theta)}}} \ ,
\la{888}
\ee
then
\be\la{3.35}
\inf_{C^{\infty}_c(\Omega)} \frac{\int_{\Omega}|\nabla u|^2dx}{\int_{\Omega}
\frac{u^2}{|x|^2}dx} \leq
\inf_{C^{\infty}_c(\ccA_{\rho,\theta})} \frac{\int_{\ccA_{\rho,\theta}}|\nabla u|^2dx}{\int_{\ccA_{\rho,\theta}}
\frac{u^2}{|x|^2}dx} <\frac{n^2}{4},
\ee
that is the Hardy inequality with constant $n^2/4$ fails in $\Omega$  if the exterior ball at zero is small enough.

{\bf \em Proof of (\ref{3.35}).} We first note that
\[
 (B_1 \setminus B_{2\rho\cos\theta}) \cap  \{ (x',x_n) : x_n<\cot\theta \; |x'| \} \subset \ccA_{\rho,\theta}.
\]
Separating variables we then conclude that
\bean
\inf_{C^{\infty}_c(\ccA_{\rho,\theta})} \frac{\int_{\ccA_{\rho,\theta}}|\nabla u|^2dx}{\int_{\ccA_{\rho,\theta}}
\frac{u^2}{|x|^2}dx} &\leq &
\inf_{f(2\rho\cos\theta)=f(1)=0}\frac{ \int_{2\rho\cos\theta}^1   f'(r)^2 r^{n-1} dr }{\int_{2\rho\cos\theta}^1   f(r)^2 r^{n-3} dr } +
\inf_{g(\theta)=g'(\pi)=0}\frac{ \int_{\theta}^{\pi}  \sin^{n-2}t  g'(t)^2 dt }{\int_{\theta}^{\pi}  \sin^{n-2}t  g(t)^2 dt }\\
&=& \Big( \frac{n-2}{2}\Big)^2 + \Big( \frac{\pi}{\ln (2\rho\cos\theta)}\Big)^2 +\lambda_1(n,\theta) 
< \frac{n^2}{4},
\eean
by assumption (\ref{888}).

\section{Improved Hardy-Sobolev inequalities for bounded domains}

In this section we shall establish improved Hardy and Hardy-Sobolev inequalities and in particular we will provide the proof of Theorem \ref{c}.
We start with the following lemma.
\begin{lemma}
Let $n\geq 3$. There exist $\sigma_n \in (0,1)$  and a constant $C_n>0$, both depending only on $n$, such that for all $\rho>0$
and all $r\leq \sigma_n\rho$ we have
\be
\int_{\cC B(\rho) \cap B(\rho e_n,r  )}|\nabla u|^2dx \geq \frac{n^2}{4}\int_{\cC B(\rho) \cap B(\rho e_n, r)}
\frac{u^2}{|x-\rho e_n|^2}dx 
+ C_n \left( \int_{\cC B(\rho) \cap B(\rho e_n, r)}  X_1^{\frac{2n-2}{n-2}}|u|^{\frac{2n}{n-2}} 
dx \right)^{\frac{n-2}{n}} , 
\la{gera:eq:1.33aaa} 
\ee
for all $u\in C^{\infty}_c(\cC B(\rho) \cap B(\rho e_n, r ))$; here $X_1=X_1(|x-\rho e_n|/3r)$.
\label{gera:thm:general_sobolev}
\end{lemma}
{\bf \em Proof.} We establish (\ref{gera:eq:1.33aaa}) for $\rho=1$, the general case will then follow by scaling.
The map
\be
S(v)= \frac{1}{|v +e_n|^2} (2v',1-|v|^2)
\la{gera:map:s}
\ee
maps conformally $\R^n_+$ onto the unit ball $B_1$. We note that
\be
|S(v)| =\frac{|v-e_n|}{|v+e_n|}.
\la{gera:S1}
\ee
Composing $S$ with the Kelvin transform $K$ we obtain that the map
\be
T(v) =(KS)(v)=  \frac{1}{ |v-e_n|^2}(2v',1-|v|^2) 
\la{gera:def_T}
\ee
maps conformally $\R^n_+$ onto $\cC B_1$.
The Jacobian determinant $JS(v)$ of $S$ can be computed explicitly and one finds
\[
|JS(v)| =\frac{2^n}{ |v +e_n|^{2n}} .
\] 
The Jacobian of the Kelvin map $K(y)$ is $|y|^{-2n}$ hence, using also (\ref{gera:S1}), the Jacobian of $T$ is
\be
|JT(v)| = | JK(Sv)|  \,\, |JS(v)| = \frac{2^n |Sv|^{-2n}}{  |v+e_n|^{2n}} = \frac{2^n}{|v-e_n|^{2n}}.
\la{gera:jt1}
\ee
Now,  simple computations give that $S^{-1}=S$ and therefore $T^{-1}=S^{-1}K^{-1}=SK$. From this we find
\[
T^{-1}(x)=\frac{1}{|x'|^2+(x_n+1)^2}(2x',|x|^2 -1)
\]  
and therefore
\be
|T^{-1}(x)| =\frac{|x-e_n|}{|x+e_n|}.
\la{gera:T1112}
\ee
Now let $r<1$ be fixed (this will be chosen later on)
and let $F\in C^{\infty}_c(T(B_r^+))$ be given. We define the function $G$ on $B_r^+$ by
\[
G(v) =F(T(v))  | (JT)(v)|^{\frac{n-2}{2n}}  = F(T(v)) \Big( \frac{2}{|v-e_n|^2} \Big) ^{\frac{n-2}{2}}.
\]
We then have by Lemma \ref{gera:lem:123},
\be
\int_{B_r^+}|\nabla G|^2  dv \geq \frac{n^2}{4} \int_{B_r^+}\frac{G^2}{|v|^2} dv  +
 \frac{1}{16r^{1/2}}\int_{B_{r}^+}\frac{G^2}{|v|^{3/2}}dv  + C_n \left( \int_{B_r+}  X_1^{\frac{2n-2}{n-2}}|G|^{\frac{2n}{n-2}}  dv \right)^{\frac{n-2}{n}}
 ,  \label{gera:josbir1}
\ee
where $X_1=X_1(|v|/r)$. We next change variables in (\ref{gera:josbir1}).

We have
\[
G(v) = 2^{\frac{n-2}{2}} F(T(v))  |v-e_n|^{2-n}
\]
and therefore
\bean
|\nabla_v G(v)|^2 & =& 2^{n-2} \Big(  |\nabla F(T(v))|^2 |v-e_n|^{2(2-n)}
+ 2 |v-e_n|^{2-n}F(T(v)) \nabla (F(T(v)) \cdot \nabla |v-e_n|^{2-n} \\
&& +F(T(v))^2 | \nabla |v-e_n|^{2-n} |^2  \Big). 
\eean
After integration over $B_r^+$ and a change of variables the first term turns out to be equal to $\int_{T(B_r^+)}|\nabla F|^2  dx$.
Integrating the other two terms yields
\begin{eqnarray*}
&&\hspace{-1.5cm}\int_{B_r^+}\Big(  
2 |v-e_n|^{2-n}F(T(v)) \nabla (F(T(v)) \cdot \nabla |v-e_n|^{2-n}
+F(T(v))^2 | \nabla |v-e_n|^{2-n} |^2  \Big)dv \\
&=&  \int_{B_r^+}\Big(  
 |v-e_n|^{2-n} \nabla (F(T(v))^2  \cdot \nabla |v-e_n|^{2-n}
+F(T(v))^2 | \nabla |v-e_n|^{2-n} |^2  \Big)dv \\
&=&  \int_{B_r^+} (F(T(v)))^2  \Big(  - \diver \big(  |v-e_n|^{2-n}  \nabla |v-e_n|^{2-n} \big)
+ \big| \nabla |v-e_n|^{2-n} \big|^2  \Big)dv \\
&=&0.
\end{eqnarray*}
We thus conclude that
\[
\int_{B_r^+}|\nabla G|^2  dv  =  \int_{T(B_r^+)}|\nabla F|^2  dx \, .
\]
Using  (\ref{gera:jt1}) and (\ref{gera:T1112}) we also find that
\[
\int_{B_r^+}\frac{G^2}{|v|^2} dv = \int_{T(B_r^+)}\frac{4F^2}{|x-e_n|^2|x+e_n|^2} dx \, .
\]
The other two integrals in (\ref{gera:josbir1}) can similarly be transformed and we conclude that (\ref{gera:josbir1})
takes the form
\bean
 \int_{T(B_r^+)}|\nabla F|^2dx &\geq& \frac{n^2}{4} \int_{T(B_r^+)}\frac{4F^2}{|x-e_n|^2|x+e_n|^2}dx   + \frac{1}{16r^{1/2}} \int_{T(B_r^+)}\frac{4F^2}{|x-e_n|^{3/2} |x+e_n|^{5/2}}dx  \\
&& +C_n \left( \int_{T(B_r+)}  X_1^{\frac{2n-2}{n-2}}|F|^{\frac{2n}{n-2}}  dv \right)^{\frac{n-2}{n}},
\eean
where $X_1=X_1( |x-e_n| /r|x+e_n|)$.
Now, it follows from (\ref{gera:T1112}) and some simple geometry that for any $r<1$
\[
T(B(r))=\{ x\in\R^n \, : \, |x-e_n| <  r|x+e_n|\} =B\big(\frac{1+r^2}{1-r^2} e_n , \frac{2r}{1-r^2}\big) \supset B(e_n, r),
\]
therefore
\be
T(B_r^+)  \supset B_1^{c}  \cap B(e_n, r)  .
\la{inclusion}
\ee
We will choose $\sigma_n \in (0,1)$ such that for all $r\leq\sigma_n$ and for all $x\in B_1^{c}\cap B(e_n, r)
\subset T(B_r^+)$ there holds
\[
\frac{n^2}{4}\frac{4}{|x-e_n|^2|x+e_n|^2} +  \frac{1}{16r^{1/2}}\frac{4}{|x-e_n|^{3/2} |x+e_n|^{5/2}}  
\geq \frac{n^2}{4|x-e_n|^2} ,
\]
or equivalently,
\be
|x-e_n|^{1/2} \geq  n^2r^{1/2}|x+e_n|^{5/2} \Big( 1-\frac{4}{|x+e_n|^2}\Big).
\la{gef}
\ee
Indeed, this is immediate for $|x+e_n| \leq 2$.  Assuming that $|x+e_n| > 2$ we 
set $t=|x-e_n|$. We then have $|x+e_n|\leq t+2$ and therefore (\ref{gef}) will follow provided
\[
n^2 r^{1/2} t^{1/2}(t+4)(t+2)^{1/2} \leq 1,
\]
for all $t\leq r$. Simple computations give that the last inequality holds true provided $ t\leq 1 / (75n^4r)$.
This will be true for all $t\leq r$ and all $r\leq\sigma_n$ if $\sigma_n$ is chosen as
\[
\sigma_n =\frac{1}{\sqrt{75}n^2}.
\]
Finally, the inequality $|x+e_n| \leq 3$ implies $X_1( |x-e_n| /r|x+e_n|)\geq X_1(|x-e_n|/3r)$. This completes the proof.
$\hfill\Box$


\noindent
{\bf \em Proof of Theorem \ref{c}.} We shall actually prove that the constant $\sigma_n$ in the statement of the theorem is the same as the constant $\sigma_n$ in the statement of Lemma~\ref{gera:thm:general_sobolev}. Without loss of generality we may assume that $\rho=1$, the general case then follows by scaling.

We first note that from Lemma \ref{gera:thm:general_sobolev} and the inclusions
\[
\Omega \cap B_{\rho} \subset \cC B(-e_n, 1) 
\cap B_{\rho}   \subset  \cC B(-e_n, 1) \cap B_{\sigma_n} .
\]
we obtain that for all $r \leq\sigma_n\rho $ there holds
\be
\int_{\xO \cap B_r} |\nabla u|^2 dx \geq   \frac{n^2}{4}  \int_{\xO \cap B_r} \frac{u^2}{|x|^2} dx 
+  C_n \bigg( \int_{\xO \cap B_r}  X_1^{\frac{2n-2}{n-2}}|u|^{\frac{2n}{n-2}}     dx \bigg)^{\frac{n-2}{n}},
\la{wasa}
\ee
for all $u \in C^{\infty}_{c}(\xO \cap B_r)$, where $X_1=X_1(|x|/3r)$.

We apply (\ref{wasa}) for $r=D$ (which is allowed since $D\leq\sigma_n\rho$) and the result follows immediately from the inclusion $\Omega\subset B_D$.
To establish the optimality of  the exponent $(2n-2)/(n-2)$,  it is enough to show the following

{\bf Claim.} If $p<(2n-2)/(n-2)$ then there is no $\sigma>1$ such that the inequality
\be
\int_{ B_1\cap B_{\rho}(e_n)} |\nabla u|^2 dx \geq   \frac{n^2}{4}  \int_{ B_1\cap B_{\rho}(e_n)} \frac{u^2}{|x-e_n|^2} dx 
+  C \left( \int_{ B_1\cap B_{\rho}(e_n)}  X_1^{p}|u|^{\frac{2n}{n-2}}     dx \right)^{\frac{n-2}{n}}  \label{gera:dbp:opt}
\ee
with $X_1=X_1(|x-e_n|/\sigma\rho)$ holds
true for some small $\rho>0$ and some $C>0$ and all $u \in C^{\infty}_{c}(\xO \cap B_{\rho})$.

Suppose to the contrary that (\ref{gera:dbp:opt}) is true for all $u\in C^{\infty}_c(B_1\cap B_{\rho}(e_n))$.
We use the conformal map $S$ defined by (\ref{gera:map:s}) to pull-back (\ref{gera:dbp:opt}) to $S^{-1}(B_1\cap B_{\rho}(e_n) )$. We write $x=Sv$ and define
\[
w(v) =u(S(v)) \Big(  \frac{2}{|v+e_n|^2}\Big)^{\frac{n-2}{2}} .
\]
Noting that $|x-e_n| =2|v|/|v+e_n|$ we obtain that
there exists $R>0$ such that the following inequality holds true for all $w \in C^{\infty}_c(B^+_R)$
\[
\int_{B^+_R} |\nabla w|^2 dv \geq   \frac{n^2}{4}  \int_{ B^+_R} \frac{w^2}{|v|^2|v+e_n|^2} dv 
+  C \left( \int_{ B^+_R}  X_1^{p}|w|^{\frac{2n}{n-2}}     dv \right)^{\frac{n-2}{n}},  
\]
or equivalently
\be
\frac{n^2}{4}\int_{B^+_R} \frac{|v|^2 +2v_n}{|v|^2|v+e_n|^2}w^2dv  +
\int_{B^+_R} |\nabla w|^2 dv \geq   \frac{n^2}{4}  \int_{ B^+_R} \frac{w^2}{|v|^2} dv 
+  C \left( \int_{ B^+_R}  X_1^{p}|w|^{\frac{2n}{n-2}}     dv \right)^{\frac{n-2}{n}},  
\la{444}
\ee
where $X_1=X_1(2|v|/ \sigma \rho |v+e_n|)$.
Now, from Lemma \ref{gera:lem:1.12im_ex} we have the inequality
\be
\int_{B_{R}^+}\frac{w^2}{|x|^{3/2}}dx \leq  4R^{1/2} \Big( 
\int_{B_{R}^+}|\nabla w|^2dx - \frac{n^2}{4}\int_{B_{R}^+}\frac{w^2}{|x|^2}dx  \Big).
\la{444a}
\ee
By taking $R$ small enough we obtain from (\ref{444}) and (\ref{444a}) that
\[
\int_{B^+_R} |\nabla w|^2 dv \geq   \frac{n^2}{4}  \int_{ B^+_R} \frac{w^2}{|v|^2} dv 
+  C \left( \int_{ B^+_R}  X_1^{p}|w|^{\frac{2n}{n-2}}     dv \right)^{\frac{n-2}{n}} , \qquad w\in C^{\infty}_c(B^+_R).
\]
This violates the optimality of the exponent $2(n-1)/(n-2)$ of Theorem~\ref{i}, concluding the proof.
$\hfill\Box$

If the  radius of the  exterior ball is small we then have

\begin{theorem}
\label{4c}
Let $n \geq 3$. There exist positive constants $\xl_n $ and $C_n$  that depend only on $n$ such that, if the radius of the exterior ball satisfies  $\rho < D / \sigma_n$ the following holds true: 
\[
\int_{\xO} |\nabla u|^2 dx + \frac{\xl_{n}}{\rho^2} \int_{\xO} u^2 dx
\geq   \frac{n^2}{4}  \int_{\xO} \frac{u^2}{|x|^2} dx \\
+  C_n \bigg( \int_{\xO}  X_1^{\frac{2n-2}{n-2}}|u|^{\frac{2n}{n-2}}     dx \bigg)^{\frac{n-2}{n}},
\]
for all $u \in C^{\infty}_{c}(\xO)$; here $X_1=X_1(|x|/ 3D)$.
If in addition $\Omega$ satisfies an interior ball condition at $0$ then the constant $n^2/4$ and
the exponent $(2n-2)/(n-2)$ of $X_1$ are sharp in both inequalities.
\end{theorem}

{\bf \em Proof.} Without loss of generality we may assume that $\rho =1$, the general case following by scaling.
We consider a  $C^{\infty}$ cutoff function $\phi(r)$ 
such that $\phi(r) = 1$ for $0\leq r \leq \sigma_n/2$ and $\phi(r) = 0$ for $ r \geq \sigma_n$, and  $\phi(|x|)$ we have $0 \leq \phi \leq 1$, $|\nabla \phi| \leq C_1$, $|\Delta \phi|< C_2$ for some constants
 depending only on $n$. We then  compute
\bean
\int_{\Omega}|\nabla u|^2 dx  &=&  \int_{\Omega}|\nabla (\phi u) +\nabla ((1-\phi)u)|^2 dx\\
&=&  \int_{\Omega}|\nabla (\phi u)|^2 dx  +\int_{\Omega} |\nabla ((1-\phi)u)|^2 dx  + 2 \int_{\Omega}\phi(1-\phi)|\nabla u|^2dx \\
&& + \int_{\Omega}(2\phi-1)\Delta\phi \,  u^2 dx \\
&\geq&  \int_{\Omega}|\nabla (\phi u)|^2 dx  +\int_{\Omega} |\nabla ((1-\phi)u)|^2 dx  - c_n\int_{\Omega}u^2 dx \\
\mbox{ (by (\ref{wasa}))} \;\; &\geq&  \frac{n^2}{4} \int_{\xO \cap B(\sigma_n)} \frac{ \phi^2 u^2}{|x|^2}dx  + 
  C_n \Big( \int_{\xO \cap B(\sigma_n)}  X_1^{\frac{2n-2}{n-2}} (|x|/3\sigma_n) |\phi u|^{\frac{2n}{n-2}}  dx \Big)^{\frac{n-2}{n}} \\
	&& + S_n \Big( \int_{\xO}  |(1-\phi) u|^{\frac{2n}{n-2}}  dx \Big)^{\frac{n-2}{n}} - c_n\int_{\Omega}u^2 dx \\
&\geq& \frac{n^2}{4} \int_{\xO} \frac{ u^2}{|x|^2}dx  + 
  C_n \Big( \int_{\xO }  X_1^{\frac{2n-2}{n-2}} (|x|/3\sigma_n) |\phi u|^{\frac{2n}{n-2}}  dx \Big)^{\frac{n-2}{n}} \\
	&& + S_n \Big( \int_{\xO} X_1^{\frac{2n-2}{n-2}} (|x|/3\sigma_n) |(1-\phi) u|^{\frac{2n}{n-2}}dx \Big)^{\frac{n-2}{n}}-c_n'\int_{\Omega}u^2 dx \\
&\geq& \frac{n^2}{4} \int_{\xO} \frac{ u^2}{|x|^2}dx  + 
  C_n'\Big( \int_{\xO}   X_1^{\frac{2n-2}{n-2}}(|x|/3D)|u|^{\frac{2n}{n-2}}  dx \Big)^{\frac{n-2}{n}} -c_n'\int_{\Omega}u^2 dx ,
\eean
where for the last inequality we used the fact that $D>\sigma_n$.

The sharpness of the constant $n^2/4$ and of the  exponent $(2n-2)/(n-2)$ follow as in the  proof of Theorems \ref{b} and \ref{c}.
\finedim


\section{Characterizing maximal potentials}
\label{section:impr}

Throughout this section we assume  that $\Omega$ satisfies both an interior and exterior ball condition at $0$.
Without loss of generality we may  assume that the exterior ball at 0 is $B(-2\rho e_n,2\rho )$ for some
$\rho>0$.

Our starting point is the following improved Hardy inequality contained in Theorem \ref{4c},
\be
 \lambda \int_{\xO} u^2 dx  + \int_{\xO} |\nabla u|^2 dx  \geq   \frac{n^2}{4}  \int_{\xO} \frac{u^2}{|x|^2} dx 
\la{spe}
\ee
for all $u \in C^{\infty}_{c}(\xO)$. We shall be interested in the problem of improvements of (\ref{spe}) and whether corresponding best constants are attained.
In connection with this we make the following definition 
\begin{definition}
\la{def:1}
A non-negative potential $V\in L^{n/2}_{\loc}(\overline{\Omega}\setminus\{0\})$ is called {\em admissible} if there exist
$\lambda\geq 0$ and $C>0$ such that 
\be
\lambda\ino u^2dx +
\ino|\nabla u|^2dx \geq  \frac{n^2}{4} \ino\frac{u^2}{|x|^2}dx +C \ino Vu^2dx, \quad u \in  C^{\infty}_{c}(\xO),
\la{in301a}
\ee
The  class of all admissible potentials for the domain $\Omega$ is denoted by $\cA(\Omega)$.
\end{definition}
For a given $V \in \cA(\Omega)$ we denote by $b(\lambda)>0$ the best constant $C$ of inequality (\ref{in301a}). We next address the question whether
there exists non-negative potentials $W \in \cA(\Omega)$ and a positive constant $C$ such that
\[
\lambda\ino u^2dx +
\ino|\nabla u|^2dx \geq  \frac{n^2}{4} \ino\frac{u^2}{|x|^2}dx +b(\lambda)\ino Vu^2dx   +C \ino Wu^2dx,          \quad u\in C^{\infty}_{c}(\xO) .
\]
In case there does not exist such a potential $W$  we say that the potential
\[
\frac{n^2}{4} \frac{1}{|x|^2} +b(\lambda)V(x) \ ,
\]
is a {\it maximal potential}. Our next goal is to characterize maximal potentials. In this direction 
for $V\in\cA(\Omega)$  and small $r>0$ we define
\be
\la{crv}
C_r(V) =\inf_{u\in C^{\infty}_{c}(\xO)}
\frac{\lambda\int_{\Omega\cap B_r}u^2dx +
\int_{\Omega\cap B_r}|\nabla u|^2dx - \frac{n^2}{4} 
\int_{\Omega\cap B_r}\frac{u^2}{|x|^2}dx}{\int_{\Omega\cap B_r}Vu^2dx}.
\ee
Since $C_r(V)$ is a non-increasing function we can define
\[
\cC^0(V)=\lim_{r\to 0+}C_r(V),
\]
which may also be equal to $+\infty$.  
This definition gives the impression that $\cC^0(V)$ might depend on the choice of $\lambda$. We will now establish that $\cC^0(V)$ is independent of $\xl$. Let us denote at the moment the infimum in   (\ref{crv}) by    $\cC^0(V, \xl)$ to express the dependence on $\xl$.
We have seen (cf. (\ref{wasa}))  that for small $r>0$ there exists a positive constant $C_n$ that depends only on $n$ so that 
 \[
\int_{\xO \cap B_r} |\nabla u|^2 dx \geq   \frac{n^2}{4}  \int_{\xO \cap B_r} \frac{u^2}{|x|^2} dx 
+  C_n \bigg( \int_{\xO \cap B_r}  X_1^{\frac{2n-2}{n-2}}|u|^{\frac{2n}{n-2}}     dx \bigg)^{\frac{n-2}{n}}, \qquad 
u\in C^{\infty}_c(\xO \cap B_r).
\]
Using Holder's inequality we conclude the existence of a positive constant $c$ independent of $r$, so that for small $r$ we have
\[
\int_{\Omega\cap B_r}|\nabla u|^2dx \geq \frac{n^2}{4} \int_{\Omega\cap B_r}\frac{u^2}{|x|^2}dx + \frac{c}{r^2}
\int_{\Omega\cap B_r}u^2 dx \; , \qquad u \in C^{\infty}_c(\Omega\cap B_r).
\]
This implies
\bean
0&\leq &\frac{ \int_{\Omega\cap B_r}|\nabla u|^2dx - \frac{n^2}{4} \int_{\Omega\cap B_r}\frac{u^2}{|x|^2}dx}{\int_{\Omega\cap B_r}Vu^2dx} \leq   \frac{ \lambda\int_{\Omega\cap B_r}u^2dx +
\int_{\Omega\cap B_r}|\nabla u|^2dx - \frac{n^2}{4} \int_{\Omega\cap B_r}\frac{u^2}{|x|^2}dx}{\int_{\Omega\cap B_r}Vu^2dx} \\
&\leq& (1+\lambda cr^2)\frac{ \int_{\Omega\cap B_r}|\nabla u|^2dx - \frac{n^2}{4} \int_{\Omega\cap B_r}\frac{u^2}{|x|^2}dx}{\int_{\Omega\cap B_r}Vu^2dx}.
\eean
Hence $C_r(V,0)\leq C_r(V, \lambda)\leq  (1+\lambda cr^2) C_r(V,0)$. Letting $r\to 0$ we conclude that $\cC^0(V)$ is indeed
independent of the choice of $\lambda \geq 0$.

\begin{definition}
We say that the potential $V\in \cA(\Omega)$ is subcritical if $\cC^0(V) =+\infty$.
\la{def:2}
\end{definition}

\begin{lemma}
\la{lem:77}
Let $V$ be a non-negative potential satisfying
\be
\ino V^{n/2}X_1^{1-n}dx <+\infty,
\la{77}
\ee
where $X_1=X_1(|x|/D)$. Then $V$ is a subcritical potential.
\end{lemma}
{\bf \em Proof.} Applying  Theorem \ref{c} we obtain that for $r>0$ small enough we have
\[
\int_{ \xO \cap B_r} |\nabla u|^2 dx \geq   \frac{n^2}{4}  \int_{\xO \cap B_{r}} \frac{u^2}{|x|^2} dx 
+  C_n \bigg( \int_{\xO \cap B_{r}}  X_1^{\frac{2n-2}{n-2}}|u|^{\frac{2n}{n-2}}     dx \bigg)^{\frac{n-2}{n}}
\]
for all $u\in C^{\infty}_c(\xO \cap B_{r})$, where $X_1=X_1(|x|/D)$. Applying H\"{o}lder inequality we then easily obtain that
\[
\frac{\int_{\Omega\cap B_r}|\nabla u|^2dx - \frac{n^2}{4} \int_{\Omega\cap B_r}\frac{u^2}{|x|^2}dx}
{\int_{\Omega\cap B_r}Vu^2dx} \geq \frac{1}{ C_n} \Big(   \int_{\Omega\cap B_r} V^{n/2}X_1^{1-n}dx \Big)^{-\frac{2}{n}}.
\]
Letting $r\to 0+$ we conclude that $\cC^0(V) =+\infty$.
 $\hfill\Box$

We shall also consider the following more general situation. Assume that $V,W_1,W_2$ are non-negative potentials in $\cA(\Omega)$ and assume that there exist
$c>0$ and a radius $R>0$ so that
\be
\int_{\Omega\cap B_{R}}W_1 u^2dx +
\int_{\Omega\cap B_{R}}|\nabla u|^2dx \geq \frac{n^2}{4} \int_{\Omega\cap B_{R}}\frac{u^2}{|x|^2}dx +
\int_{\Omega\cap B_{R}}W_2 u^2dx + c\int_{\Omega\cap B_{R}}V u^2dx 
\la{w1w2}
\ee
for all $u\in C^{\infty}_c(\Omega\cap B_{R})$. For $0<r<R$ we define
\[
C_r(W_1,W_2;V) =\inf_{C^{\infty}_c(\Omega\cap B_{r})}\frac{\int_{\Omega\cap B_{r}}W_1 u^2dx +\int_{\Omega\cap B_{r}}|\nabla u|^2dx - \frac{n^2}{4} \int_{\Omega\cap B_{r}}\frac{u^2}{|x|^2}dx -
\int_{\Omega\cap B_{r}}W_2 u^2dx}{\int_{\Omega\cap B_{r}}V u^2dx}
\]
and we denote
\[
\cC^0(W_1,W_2;V) =\lim_{r\to 0} C_r(W_1,W_2;V).
\]

We next show that subcritical potentials do not affect the concentration level $\cC^0(V)$. More precisely we have
\begin{lemma}
Let $V,W_1,W_2$ be non-negative potentials in $\cA(\Omega)$ and assume that for some $R>0$ there exists $c>0$ such that (\ref{w1w2}) holds true.
If in addition $W_1,W_2$ are subcritical then $\cC^0(W_1,W_2;V) =\cC^0(V)$.
\label{lem:w1w2}
\end{lemma}
{\bf \em Proof.} The subcriticality of $W_i$ implies that for small $r>0$ there holds
\be
\la{ww1}
\int_{\Omega\cap B_{r}}|\nabla u|^2dx \geq \frac{n^2}{4} \int_{\Omega\cap B_{r}}\frac{u^2}{|x|^2}dx +
c_r(W_i)\int_{\Omega\cap B_{r}}W_i u^2dx \; , \quad u \in C^{\infty}_c(\Omega\cap B_{r}) \ ,
\ee
with $\lim_{r\to 0} c_r(W_i) =+\infty$. From inequalities (\ref{w1w2}) and (\ref{ww1}) follows that
for $r>0$ small enough we have
\bean
&&\hspace{-4cm}\Big(1 -\frac{1}{c_r(W_2)}\Big)\frac{\int_{\Omega\cap B_{r}}|\nabla u|^2dx - \frac{n^2}{4} \int_{\Omega\cap B_{r}}\frac{u^2}{|x|^2}dx}{ \int_{\Omega\cap B_{r}}V u^2dx}\\
&\leq&
\frac{\int_{\Omega\cap B_{r}}|\nabla u|^2dx - \frac{n^2}{4} \int_{\Omega\cap B_{r}}\frac{u^2}{|x|^2}dx +
 \int_{\Omega\cap B_{r}}(W_1-W_2) u^2dx}{ \int_{\Omega\cap B_{r}}V u^2dx} \\
&\leq&\Big(1 +\frac{1}{c_r(W_1)}\Big)
\frac{\int_{\Omega\cap B_{r}}|\nabla u|^2dx - \frac{n^2}{4} \int_{\Omega\cap B_{r}}\frac{u^2}{|x|^2}dx }{ \int_{\Omega\cap B_{r}}V u^2dx} \ .
\eean
This implies that 
\[
\Big(1 -\frac{1}{c_r(W_2)}\Big) c_r(V) \leq   C_r(W_1,W_2;V) \leq \Big(1 +\frac{1}{c_r(W_1)}\Big) c_r(V),
\]
and the result follows by letting $r\to 0+$.

 $\hfill\Box$

Given $u\in C^{\infty}_c(\Omega)$ we define the function $w$ by
\be
u(x) =  \frac{|x+\rho e_n|^2- \rho^2}{|x|^{\frac{n}{2}}   |x+2\rho e_n|^{\frac{n}{2}}} w(x).
\la{uw}
\ee
Then $w\in C^{\infty}_c(\Omega)$ by our assumption that the exterior ball at zero is $B(-2\rho e_n,2\rho)$.
After some computations and using integration by parts we arrive at
\[
\ino|\nabla u|^2dx =\ino  \frac{\big(  |x+\rho e_n|^2- \rho^2 \big)^2 }{|x|^{n}   |x+2\rho e_n|^{n}}|\nabla w|^2dx +
 n^2 \rho^2 \ino  \frac{\big(  |x+\rho e_n|^2- \rho^2 \big)^2 }{|x|^{n+2}   |x+2\rho e_n|^{n+2}} w^2dx ,
\]
so inequality (\ref{in301a}) is written
 \bean
\lambda\ino \frac{\big( |x+ \rho e_n|^2- \rho^2 \big)^2 }{|x|^{n}   |x+2\rho e_n|^{n}}   w^2dx &+&
\ino  \frac{\big(  |x+\rho e_n|^2- \rho^2 \big)^2 }{|x|^{n}   |x+2\rho e_n|^{n}}|\nabla w|^2dx
+  n^2 \rho^2 \ino  \frac{\big(  |x+\rho e_n|^2-\rho^2 \big)^2 }{|x|^{n+2}   |x+2\rho e_n|^{n+2}} w^2dx \\
& \geq & \frac{n^2}{4}\ino  \frac{\big(  |x+\rho e_n|^2-\rho^2 \big)^2 }{|x|^{n+2}   |x+2\rho e_n|^{n}} w^2dx 
+ c \ino \frac{\big( |x+\rho e_n|^2-\rho^2 \big)^2 }{|x|^{n}   |x+2\rho e_n|^{n}}  V w^2 dx
\eean
which can also take the equivalent form
 \bea
\lambda\ino \frac{\big( |x+\rho e_n|^2- \rho^2 \big)^2 }{|x|^{n}   |x+2\rho e_n|^{n}}   w^2dx &+&
\ino  \frac{\big(  |x+\rho e_n|^2-\rho^2 \big)^2 }{|x|^{n}   |x+2\rho e_n|^{n}}|\nabla w|^2dx \nn \\
&  &\hspace{-3cm}\geq \frac{n^2}{4}\ino  \frac{\big(  |x+ \rho e_n|^2-\rho^2 \big)^2 (|x|^2+4\rho x_n)}{|x|^{n+2}   |x+2\rho e_n|^{n+2}} w^2dx 
+ c \ino \frac{\big( |x+\rho e_n|^2-\rho^2 \big)^2 }{|x|^{n}   |x+2\rho e_n|^{n}}  V w^2 dx.
\label{w}
\eea
It is clear from the above that (\ref{w}) is valid for all functions $w\in C^{\infty}_c(\Omega)$ and moreover,
for a fixed $\lambda \geq 0$ the best constants $c$ for inequalities (\ref{in301a}) and (\ref{w}) coincide.
That common best constant shall be denoted by $b=b(\lambda)$.

Denoting
\[
\phi(x)=\frac{|x+ \rho e_n|^2-\rho^2 }{|x|^{\frac{n}{2}}   |x+2\rho e_n|^{\frac{n}{2}}} 
\]
it then follows that (cf. (\ref{crv}))
\[
 C_r(V) = 
\inf_{ W^{1,2}_0( \Omega\cap B_r ; \phi^2)}
\hspace{-0.3cm}\frac{\lambda \int_{\Omega\cap B_r}\phi^2   w^2dx +
\int_{\Omega\cap B_r} \phi^2|\nabla w|^2dx -
 \frac{n^2}{4}\int_{\Omega\cap B_r} \frac{|x|^2+4\rho x_n}{|x|^{2}   |x+2\rho e_n|^{2}} \phi^2 w^2dx}
{  \int_{\Omega\cap B_r} V \phi^2  w^2 dx}
\]
where $W^{1,2}_0( \Omega\cap B_r ; \phi^2)$ denotes the closure of $C^{\infty}_c(\Omega\cap B_r)$ under the norm
\[
\int_{\Omega\cap B_r}\phi^2|\nabla w|^2dx +
 \int_{\Omega\cap B_r} \phi^2  w^2dx .
\]

We shall now see a simpler way for expressing $\cC^0(V)=\lim_{r\to 0+}C_r(V)$ in terms of the weight $\phi^2$. For this we define
\[
C_r(V;\phi^2) =\inf_{W^{1,2}_0( \Omega\cap B_r ; \phi^2)}
\frac{ \int_{\Omega\cap B_r}  \phi^2|\nabla w|^2dx}{ \int_{\Omega\cap B_r} V  \phi^2  w^2dx}.
\]
and
\[
\cC^0(V;\phi^2) =\lim_{r\to 0+} C_r(V;\phi^2).
\]

\begin{lemma}
Let $V$ be a non-negative potential in $\cA(\Omega)$. Then
$\cC^0(V) =\cC^0(V;\phi^2)$.
\label{lem:fr_ir}
\end{lemma}
{\bf \em Proof.} For the sake of simplicity we assume that $\rho=1$; the general case then follows by scaling.
On the one hand we have for small $r>0$ that
\[
cr^{-2}\int_{\Omega\cap B_r}\phi^2w^2 dx \leq  \int_{\Omega\cap B_r} \phi^2|\nabla w|^2dx -
 \frac{n^2}{4}\int_{\Omega\cap B_r} \frac{|x|^2+4x_n}{|x|^{2}   |x+2e_n|^{2}} \phi^2 w^2dx \; , \qquad
u\in C^{\infty}_c(\Omega\cap B_r)
\]
for some universal constant $c>0$, which implies the inequality
\be
\frac{\lambda \int_{\Omega\cap B_r}\phi^2   w^2dx +
\int_{\Omega\cap B_r} \phi^2|\nabla w|^2dx -
 \frac{n^2}{4}\int_{\Omega\cap B_r} \frac{|x|^2+4x_n}{|x|^{2}   |x+2e_n|^{2}} \phi^2 w^2dx}
{  \int_{\Omega\cap B_r} \phi^2 V w^2 dx}
\leq
\Big( 1+\frac{\lambda}{c}r^2\Big)\frac{\int_{\Omega\cap B_r} \phi^2|\nabla w|^2dx }{  \int_{\Omega\cap B_r} \phi^2 V w^2 dx}.
\la{ineq1}
\ee
On the other hand, since
\[
0\leq \frac{|x|^2 +4x_n}{|x|^2 |x+2e_n|^2} \leq \frac{c}{|x|},
\]
we have the inequality
\[
cr^{-1}\int_{\Omega\cap B_r} \frac{|x|^2+4x_n}{|x|^{2}   |x+2e_n|^{2}} \phi^2 w^2dx \leq \int_{\Omega\cap B_r} \phi^2|\nabla w|^2dx \, , \qquad  w\in C^{\infty}_c(\Omega\cap B_r) 
\]
which in turn implies
\be
\frac{\int_{\Omega\cap B_r} \phi^2|\nabla w|^2dx }{  \int_{\Omega\cap B_r} \phi^2 V w^2 dx}
(1-cr)
\leq 
\frac{\lambda \int_{\Omega\cap B_r}\phi^2   w^2dx +
\int_{\Omega\cap B_r} \phi^2|\nabla w|^2dx -
 \frac{n^2}{4}\int_{\Omega\cap B_r} \frac{|x|^2+4x_n}{|x|^{2}   |x+2e_n|^{2}} \phi^2 w^2dx}
{  \int_{\Omega\cap B_r} \phi^2 V w^2 dx} 
\la{ineq2}
\ee
The result follows by combining inequalities (\ref{ineq1}) and (\ref{ineq2}) 
and letting $r\to 0+$. $\hfill\Box$

One important consequence of subcriticality is the following compactness property.
\begin{lemma}
Assume that the positive potential $V\in\cA(\Omega)$ is subcritical. Then for any sequence $(w_k)$ which is bounded in $W^{1,2}_0(\Omega ;\phi^2)$ there exists a subsequence, also denoted by $(w_k)$, and a function $w_0\in W^{1,2}_0(\Omega ;\phi^2)$ so that
\bean
\ia && w_k \rightharpoonup w_0 \mbox{ in } W^{1,2}_0(\Omega ;\phi^2)\\
\ib && \int_{\Omega}\phi^2V(w_k-w_0)^2dx \longrightarrow 0 .
\eean
\label{compactness}
\end{lemma}
{\bf \em Proof.} Part (i) is standard. To prove (ii) we may assume without loss of generality that $w_0=0$.
We consider a small $r>0$ and a smooth cut-off function $\psi$ such that $\psi=1$ in $B_{r/2}$ and $\psi=0$ outside $B_r$.
We then have
\bean
 \int_{\Omega}\phi^2|\nabla w_k|^2 dx &=& \int_{\Omega}\phi^2|\nabla (\psi w_k)  +  \nabla ((1-\psi) w_k)  |^2 dx \\
&=&
 \int_{\Omega}\phi^2 |\nabla (\psi w_k)|^2dx +  \int_{\Omega}\phi^2 |\nabla \big( (1-\psi) w_k\big)|^2dx  +
2 \int_{\Omega}\phi^2 \psi(1-\psi)|\nabla w_k|^2dx \\
&&  + \int_{\Omega}\phi^2 (1-2\psi)w_k\nabla\psi\cdot \nabla w_k dx - \int_{\Omega}\phi^2  |\nabla \psi|^2 w_k^2dx \\
&\geq& \int_{\Omega}\phi^2 |\nabla (\psi w_k)|^2dx +o(1) \\
&\geq& C_r(V;\phi^2) \int_{\Omega}\phi^2 V \psi^2w_k^2dx +o(1)\\
&=&  C_r(V;\phi^2) \int_{\Omega}\phi^2 V w_k^2dx +o(1) , 
\eean
that is
\be
\int_{\Omega} V\phi^2  w_k^2dx \leq \frac{1}{C_r(V;\phi^2)} \int_{\Omega}\phi^2|\nabla w_k|^2 dx +o(1) \ .
\la{pas}
\ee
The result follows by noting that the RHS of (\ref{pas}) can be made arbitrarily small by choosing $r>0$ small enough.
$\hfill\Box$

We can now state and prove the main result of this section.
\begin{theorem}
Let $V\in\cA(\Omega)$ and $\lambda\geq 0$ be given and let $b(\lambda)$ be the best constant for the inequality
\be
\lambda\ino \phi^2 w^2dx +
\ino \phi^2|\nabla w|^2dx \geq  \frac{n^2}{4} \ino   \frac{|x|^2+4\rho x_n}{|x|^{2} |x+2\rho e_n|^{2}}  \phi^2 w^2
dx +b(\lambda)\ino V \phi^2 w^2dx , \qquad w\in W^{1,2}_0(\Omega;\phi^2).
\la{eq:w:no1}
\ee
If in addition
\[
b(\lambda)<\cC^0(V)
\]
then the best constant $b(\lambda)$ in (\ref{eq:w:no1}) is realized by a function $w_0\in W^{1,2}_0(\Omega;\phi^2)$. In particular the best constant $b(\lambda)$ is realized if the potential $V$ is subcritical.
\label{gera:thm:exmin}
\end{theorem}
{\bf \em Proof.} We denote
\[
Q(x)=\frac{n^2}{4}  \frac{|x|^2+4\rho x_n}{|x|^{2} |x+2\rho e_n|^{2}}
\]
Then it is easily seen that
\[
0\leq Q\leq \frac{c}{|x|}, \qquad x\in\Omega,
\]
which implies that $Q$ is a subcritical potential by Lemma \ref{lem:77}.
We consider a minimizing sequence $(w_k)$ for (\ref{eq:w:no1}) and without loss of generality we assume that
\be
\ino Q \phi^2  w_k^2dx  +  b(\lambda) \ino V \phi^2  w_k^2 dx
=1, \qquad   \lambda\ino \phi^2  w_k^2dx + \ino  \phi^2 |\nabla w_k|^2dx  \longrightarrow 1, \mbox{ as }k\to \infty.
\label{g11}
\ee
Since $(w_k)$ is bounded in $W^{1,2}_0(\Omega ; \phi^2)$, it has a subsequence, which we assume is $(w_k)$ itself, which converges weakly
to some $w_0\in W^{1,2}_0(\Omega ; \phi^2)$. We define $v_k=w_k-w_0$. 

We consider a small enough $r>0$ so that $C_r(V;\phi^2)>b(\lambda)$ and a smooth cut-off function $\psi$
such that $\psi=1$ in $B_{r/2}$ and $\psi=0$ outside $B_r$. Arguing as in the proof of Lemma \ref{compactness} we have
\bea
\int_{\Omega}\phi^2 |\nabla v_k|^2 dx &=&  \int_{\Omega}\phi^2 |\nabla [\psi v_k +(1-\psi)v_k] |^2dx \nn \\
&\geq&  \int_{\Omega}\phi^2 |\nabla (\psi v_k)|^2dx  +o(1) \nn \\
&\geq&  C_r(V;\phi^2) \int_{\Omega}\phi^2 V \psi^2 v_k^2 dx +o(1) \nn  \\
&=&  C_r(V;\phi^2) \int_{\Omega}\phi^2 V v_k^2 dx +o(1) \ .  \la{equ}
\eea
Now, substituting $w_k=v_k+w_0$ in the normalization relations (\ref{g11}) and using Lemma \ref{compactness} we obtain
\be
\ino \phi^2 Q w_0^2dx  +  b(\lambda) \ino \phi^2 V w_0^2 dx + b(\lambda) \ino \phi^2 V v_k^2 dx  =1 +o(1)
\la{norm1}
\ee
and
\be
\lambda\ino \phi^2  w_0^2dx + \ino  \phi^2 |\nabla w_0|^2dx + \ino  \phi^2 |\nabla v_k|^2dx = 1+o(1).
\la{norm2}
\ee
From (\ref{equ}) and (\ref{norm1}) we obtain
\be
\int_{\Omega}\phi^2 |\nabla v_k|^2 dx \geq \frac{ C_r(V;\phi^2)}{b(\lambda)}\Big(  1-    \ino \phi^2 Q w_0^2dx  +  b(\lambda) \ino V \phi^2  w_0^2 dx        \Big) +o(1).
\la{norm7}
\ee
Moreover using (\ref{eq:w:no1}) for $w=w_0$ we obtain from (\ref{norm2}) that
\be
\ino  \phi^2 |\nabla v_k|^2dx \leq 1 -\int_{\Omega}\phi^2Qw_0^2dx -b(\lambda)\ino \phi^2Vw_0^2dx +o(1) .
\la{norm3}
\ee
From (\ref{norm7}) and (\ref{norm3}) we conclude that
\[
\Big(  1- \frac{ C_r(V;\phi^2)}{b(\lambda)}\Big) \bigg(      1-    \ino \phi^2 Q w_0^2dx  +  b(\lambda) \ino V \phi^2  w_0^2 dx       \bigg) \geq 0.
\]
Since $C_r(V;\phi^2)>b(\lambda)$, this implies that
\[
 \ino \phi^2 Q w_0^2dx  +  b(\lambda) \ino \phi^2 V w_0^2 dx \geq 1.
\]
But by lower semicontinuity,
\[
 \ino \phi^2 Q w_0^2dx  +  b(\lambda) \ino \phi^2 V w_0^2 dx \leq 1.
\]
Hence $w_0$ is a minimizer. $\hfill\Box$

The next theorem is an immediate consequence of Theorem \ref{gera:thm:exmin}
\begin{theorem}
\label{d}
Let $V$ be a non-negative potential in $\cA(\Omega)$. (a) Let $\lambda\geq 0$ and $b(\lambda)>0$ be such that
\be
\lambda\ino u^2dx +
\ino|\nabla u|^2dx \geq  \frac{n^2}{4} \ino\frac{u^2}{|x|^2}dx +b(\lambda)\ino Vu^2dx, \quad u\in C_{c}^{\infty}(\Omega),
\la{in30}
\ee
where $b(\lambda)$ is the best constant. If in addition
\[
b(\lambda)<\cC^0(V) \ ,
\]
then the potential $n^2/4|x|^2 +b(\lambda)V(x)$ is a maximal potential, that is inequality
(\ref{in30}) cannot be improved by adding a non-negative potential $W$ in the RHS. \nl
(b) If $V$ is a subcritical potential then there exist $\lambda\geq 0$ and a best constant $b(\lambda)>0$ such that (\ref{in30}) is true. Moreover the potential $n^2/4|x|^2 +b(\lambda)V(x)$ is a maximal potential.
\end{theorem}
{\bf \em Proof.} $(a)$ Suppose that  $n^2/4|x|^2 +b(\lambda)V(x)$ is not a maximal potential, that is there exists a non-trivial potential
$W\geq 0$ in $\cA(\Omega)$ such that
\[
\lambda\ino u^2dx +
\ino|\nabla u|^2dx \geq  \frac{n^2}{4} \ino\frac{u^2}{|x|^2}dx +b(\lambda)\ino Vu^2dx + \ino Wu^2dx
\]
holds true for all $u\in H^1_0(\Omega)$. Using the transformation (\ref{uw}) this is equivalently written as
\[
\lambda\ino \phi^2 w^2dx +
\ino \phi^2|\nabla w|^2dx \geq   \ino Q\phi^2  w^2
dx +b(\lambda)\ino \phi^2Vw^2dx  +\ino W\phi^2w^2 dx , \qquad w\in W^{1,2}_0(\Omega;\phi^2).
\]
Using $w=w_0$ where $w_0$ is the minimizer from Theorem \ref{gera:thm:exmin} we conclude that $\ino W\phi^2w_0^2 dx\leq 0$, which is a contradiction. \nl
Part $(b)$ is an immediate consequence of part $(a)$ since any subcritical potential $V$ is in $\cA(\Omega)$ and satisfies
$\cC^0(V) =+\infty$.
$\hfill\Box$


\section{Logarithmic improvements and maximal potentials}
\la{sub:log:impr}

Throughout this section we continue to assume  that $\Omega$ satisfies both an interior and exterior ball condition at $0$.
We also continue to assume that the exterior ball at 0 is $B(-2\rho e_n,2\rho)$ for some $\rho>0$.

In this section we will provide the proofs of Theorems \ref{f}, \ref{g} and also study maximal potentials in the context of logarithmic improvements of Hardy
inequality.

\subsection{Logarithmic improvements}

To prove Theorems  \ref{f} and \ref{g} we first establish the following lemmas:
\begin{lemma}
Let $n\geq 2$. There exists a positive constant  $\sigma_n$ depending only on $n$ such that for all $\rho>0$ and all $r\leq\sigma_n \rho$ we have
\bea
\ia && \int_{\cC B(\rho) \cap B(\rho e_n, r)} |\nabla u|^2 dx \geq
\frac{n^2}{4} \int_{\cC B(\rho) \cap B(\rho e_n,r)} \frac{u^2}{|x-\rho e_n|^2}dx  \nonumber\\[0.2cm]
&& +
\frac{1}{4}\sum_{i=1}^{\infty} \int_{\cC B(\rho) \cap B(\rho e_n, r)} \frac{u^2}{|x-\rho e_n|^2} X_1^2\ldots X_i^2 dx \nonumber\\[0.2cm]
\ib && \mbox{If in addition $n\geq 3$ there exists a constant $C_n$ depending only on $n$ such that for all $m \in \N$} \nonumber \\[0.2cm]
&& \int_{\cC B(\rho) \cap B(\rho e_n,  r)}|\nabla u|^2dx  \geq\frac{n^2}{4}\int_{\cC B(\rho) \cap B(\rho e_n, r)}
\frac{u^2}{|x-\rho e_n|^2}dx  \nonumber \\
&& \hspace{3cm} 
+ \frac{1}{4}\sum_{i=1}^m \int_{\cC B(\rho) \cap B(\rho e_n, r)}\frac{u^2}{|x-\rho e_n|^2}X_1^2\ldots X_i^2 dx  \nonumber \\  
&& \hspace{3cm} + C_n \left( \int_{\cC B(\rho) \cap B(\rho e_n, r)} 
 (X_1\ldots X_{m+1})^{\frac{2n-2}{n-2}}|u|^{\frac{2n}{n-2}}     dx \right)^{\frac{n-2}{n}}.
\la{gera:eq:1.33}
\eea
Both inequalities are valid for all $u\in C^{\infty}_c(\cC B(\rho) \cap B(\rho e_n, r))$ and
in both cases $X_i=X_i(|x-\rho e_n|/ (3\kappa r))$. 
\label{gera:general_k_sobolev}
\end{lemma}
{\bf \em Proof.} To prove (i) it is enough to consider the case $\rho=1$.
 We fix $r<1$ and we apply Lemma \ref{gera:lem:1.12im_ex}. Changing variables via $T$ (cf. (\ref{gera:def_T})) we obtain
\bea
&&\int_{T(B_r^+)}|\nabla u|^2dx \geq \frac{n^2}{4} \int_{T(B_r^+)}\frac{4u^2}{|x-e_n|^2|x+e_n|^2}dx  \nonumber \\
&&+\frac{1}{4}\sum_{i=1}^{\infty} \int_{T(B_r^+)} \frac{4u^2\, X_1^2 \ldots X_i^2}{|x-e_n|^2|x+e_n|^2}   dx 
 + \frac{1}{8r^{1/2}} \int_{T(B_r^+)}\frac{4u^2}{|x-e_n|^{3/2} |x+e_n|^{5/2}}dx
\la{33333}
\eea
for all $u\in C^{\infty}_c(T(B_r^+))$; here $X_k=X_k( |x-e_n| / (\kappa r|x+e_n|)$. As already noted (cf. (\ref{inclusion})) we have 
\[
\cC B_1  \cap B(e_n, r) \subset T(B_r^+) .
\]
so integrals in (\ref{33333}) can be taken over $\cC B_1  \cap B(e_n, r)$.

Once again it is enough to find $\sigma_n<1$ such that for all $r\leq\sigma_n$ 
and all $x\in \cC B_1\cap B(e_n, r)$ there holds
\bean
&&\frac{n^2}{4}\frac{4}{|x-e_n|^2|x+e_n|^2} +\frac{1}{4}\sum_{i=1}^{\infty} \frac{4}{|x-e_n|^2|x+e_n|^2}  X_1^2 \ldots X_i^2  \\
&& +  \frac{4}{8r^{1/2}|x-e_n|^{3/2} |x+e_n|^{5/2}}  
\geq \frac{n^2}{4|x-e_n|^2} +\frac{1}{4|x-e_n|^2}\sum_{i=1}^{\infty} X_1^2 \ldots X_i^2. 
\eean
or equivalently,
\[
|x-e_n|^{1/2} \geq \frac{1}{2} r^{1/2} |x+e_n|^{5/2} \Big( 1-\frac{4}{|x+e_n|^2} \Big)
\Big(n^2 +\sum_{i=1}^{\infty} X_1^2 \ldots X_i^2 \Big).
\]
This is immediate if  $|x+e_n| \leq 2$. Assuming $|x+e_n| \geq 2$ we will actually establish the stronger inequality
(recalling that $\sum_k X_1^2\ldots X_k^2 \leq 1/4$)
\be
|x-e_n|^{1/2} \geq \frac{n^2 +\frac{1}{4}}{2} r^{1/2} |x+e_n|^{5/2} \Big( 1-\frac{4}{|x+e_n|^2} \Big).
\la{gera:1212}
\ee
But this is almost the same as inequality (\ref{gef}), the only difference being that in the place of $n^2$ we now have $(n^2 +1/4)/2$; we omit further details.  

The proof of (ii) is analogous, but we now use Lemma \ref{gera:lem_mlogs}
instead of Lemma \ref{gera:lem:1.12im_ex}. Again, we may take $\rho=1$. We then fix $r<1$ and changing variables via $T$ we obtain
\bea
&& \int_{T(B_r^+)}|\nabla u|^2dx \geq \frac{n^2}{4} \int_{T(B_r^+)}\frac{4u^2}{|x-e_n|^2|x+e_n|^2}dx  
+\frac{1}{4}\sum_{i=1}^{m} \int_{T(B_r^+)} \frac{4u^2\, X_1^2 \ldots X_i^2}{|x-e_n|^2|x+e_n|^2}   dx  \nonumber  \\
&& + \frac{1}{16r^{1/2}} \int_{T(B_r^+)}\frac{4u^2}{|x-e_n|^{3/2} |x+e_n|^{5/2}}dx  +  C_n \left( \int_{T(B_r^+)} 
 (X_1\ldots X_{m+1})^{\frac{2n-2}{n-2}}|u|^{\frac{2n}{n-2}}   dx \right)^{\frac{n-2}{n}}
\la{gera:cy33}
\eea
for all $u\in C^{\infty}_c(T(B_r^+))$; here $X_k=X_k( |x-e_n| / (\kappa r|x+e_n|)$. As in the proof of pert (i), this is also true if the integrals are taken over $\cC B_1  \cap B(e_n, r)$ and $u\in C^{\infty}_c(\cC B_1  \cap B(e_n, r)) $

Hence the result will follow once we establish for all $x\in \cC B_1\cap B(e_n, r)$the inequality
\bean
&&\frac{n^2}{4}\frac{4}{|x-e_n|^2|x+e_n|^2} +\frac{1}{4}\sum_{i=1}^{m} \frac{4}{|x-e_n|^2|x+e_n|^2}  X_1^2 \ldots X_i^2  \\
&& +  \frac{4}{16r^{1/2}|x-e_n|^{3/2} |x+e_n|^{5/2}}  
\geq \frac{n^2}{4|x-e_n|^2} +\frac{1}{4|x-e_n|^2}\sum_{i=1}^{m} X_1^2 \ldots X_i^2. 
\eean
This is equivalent to
\[
|x-e_n|^{1/2} \geq  r^{1/2} |x+e_n|^{5/2} \Big( 1-\frac{4}{|x+e_n|^2} \Big)
\Big(n^2 +\sum_{i=1}^{m} X_1^2 \ldots X_i^2 \Big).
\]
The argument now goes as in part (i); we omit further details. $\hfill\Box$

\noindent
{\bf\em Proof of Theorem \ref{f}}.  Without loss of generality we assume that $\rho=1$.
We use  part (i) of Lemma \ref{gera:general_k_sobolev} for $r=\sigma_n$ making a translation of
(\ref{gera:eq:1.33}) by $- e_n$. We obtain
\[
\int_{\cC B(- e_n,1) \cap B(\sigma_n)} |\nabla u|^2 dx \geq
\frac{n^2}{4} \int_{\cC B(- e_n,1) \cap B(\sigma_n)} \frac{u^2}{|x|^2}dx  
 + \frac{1}{4}\sum_{i=1}^{\infty} \int_{\cC B(- e_n, 1) \cap B(\sigma_n)} \frac{u^2}{|x|^2} X_1^2\ldots X_i^2 dx ,
\]
for all $u\in C^{\infty}_c(\cC B(- e_n, 1) \cap B(\sigma_n))$, where $X_i=X_i(|x|/(3\kappa \sigma_n))$.
Since $\Omega \subset \cC B(- e_n, 1) \cap B(\sigma_n)$, the result follows. 
\finedim

\begin{theorem}
\label{6f}

Let $n \geq 2$. There exists a positive constant $\xl_n $  that depends only on $n$ such that, if the radius of the exterior ball satisfies  $\rho < D / \sigma_n$, with $\xs_n$ as in Theorem \ref{g}, the following holds true: 
\bean
 \int_{\xO} |\nabla u|^2 dx + \frac{\xl_{n}}{\rho^2} \int_{\xO} u^2 dx 
&\geq&   \frac{n^2}{4}  \int_{\xO} \frac{u^2}{|x|^2} dx + \frac{1}{4} \sum_{i=1}^{\infty} \int_{\xO} \frac{u^2}{|x|^2}X_1^2\ldots X_i^2 dx  
\eean
for all $u \in C^{\infty}_{c}(\xO)$; here $X_i=X_i(|x|/ (3 \xk D))$.
If in addition $\Omega$ satisfies an interior ball condition at $0$ then the constants $1/4$ are sharp  at each step.

\end{theorem}

{\bf\em Proof.}
We argue as in the proof of Theorem \ref{c}. This time however we also  use the global
estimate $\sum_i{X_1^2\ldots X_i^2}\leq 1/4$ in order to estimate uniformly the constant in front of the $L^2$ term. We omit further details.
The sharpness of the constants $1/4$ has already been proved. $\hfill\Box$

\noindent
{\bf\em Proof of Theorem \ref{g}} We argue as in the proof of Theorem \ref{f}, using now part (ii) of Lemma \ref{gera:general_k_sobolev}.
To prove the sharpness of the constants $1/4$ and the exponent $(2n-2)/(n-2)$ we argue as in the proof of Theorem~\ref{c};
we omit the details. $\hfill\Box$

In  case the exterior ball is small, working as in Theorem \ref{4c} we have the following

\begin{theorem}
\label{6g}
Let $n \geq 3$. There exist positive constants $\xl_n $ and $C_n$  that depend only on $n$ such that, if the radius of the exterior ball satisfies  $\rho < D / \sigma_n$, with $\xs_n$ as in Theorem \ref{g},  then  for any $m\in\N$  the following holds true: 
\[
 \int_{\xO} |\nabla u|^2 dx + \frac{\xl_{n}}{\rho^2} \int_{\xO} u^2 dx
\geq   \frac{n^2}{4}  \int_{\xO} \frac{u^2}{|x|^2} dx + \frac{1}{4} \sum_{i=1}^m \int_{\xO} \frac{u^2}{|x|^2}X_1^2\ldots X_i^2 dx 
+  C_n \left( \int_{\xO}  (X_1\ldots X_{m+1})^{\frac{2n-2}{n-2}}|u|^{\frac{2n}{n-2}}     dx \right)^{\frac{n-2}{n}} ,
\]
for all $u \in C^{\infty}_{c}(\xO)$; here  $X_i=X_i(|x|/ (3 \xk D))$.
If in addition $\Omega$ satisfies an interior ball condition at $0$ then the exponents $(2n-2)/(n-2)$ of $X_i$ are also sharp.
\end{theorem}

\subsection{Maximal logarithmic potentials}

Here we characterize maximal potentials in the context of logarithmically improved Hardy inequalities.

Our starting point in this subsection is the improved Hardy inequality contained in Theorem \ref{g},
\be\la{62.13}
 \lambda \int_{\xO} u^2 dx  + \int_{\xO} |\nabla u|^2 dx 
\geq   \frac{n^2}{4}  \int_{\xO} \frac{u^2}{|x|^2} dx + \frac{1}{4} \sum_{i=1}^m \int_{\xO} \frac{u^2}{|x|^2}X_1^2\ldots X_i^2 dx 
\ee
for all $u \in C^{\infty}_{c}(\xO)$; here $X_i=X_i(|x|/3 \xk \tilde D)$ where $\tilde{D} \geq D$. 
We shall be interested in the problem of improvements of (\ref{62.13}) and whether the corresponding best constants are attained.

The analysis that will follow is analogous to that of Section \ref{section:impr}; for this reason we shall avoid the details in cases where the arguments are quite similar.

\begin{definition}
A non-negative potential $V\in L^{n/2}_{\loc}(\overline{\Omega}\setminus\{0\})$ is called $m$-{\em admissible} if there exist
$\lambda\geq 0$, $\tilde D\geq D$ and $C>0$ such that
\be
\lambda\ino u^2dx +
\ino|\nabla u|^2dx \geq  \frac{n^2}{4} \ino\frac{u^2}{|x|^2}dx 
+ \frac{1}{4} \sum_{i=1}^m \int_{\xO} \frac{u^2}{|x|^2}X_1^2\ldots X_i^2 dx+
C\ino Vu^2dx, \quad u \in C^{\infty}_{c}(\xO) \ ,
\la{inma}
\ee
where $X_i=X_i(|x|/3 \xk \tilde{D})$. The class of all $m$-admissible potentials for the domain $\Omega$ is denoted  by $\cA_m(\Omega)$.
\end{definition}
We note that there is a big variety of $m$-admissible potentials. For example if $V$ satisfies 
\be
\ino V^{\frac{n}{2}}(X_1 \ldots X_{m+1})^{1-n}dx <+\infty,
\la{inexa}
\ee
where $X_i=X_i(|x|/3 \xk D)$, $i=1,\ldots, m+1$, then $V$ is $m$-admissible by
Theorem~\ref{g}.

For a given $V \in \cA_m(\Omega)$ we denote by $b_m(\lambda)>0$ the best constant $C$ of inequality (\ref{inma}). We next address the question whether inequality (\ref{inma}) with best constant $b_m(\xl)$ can be further improved. That is, whether there exists
potential $W \in \cA_m(\Omega)$ and a positive constant $C$ such that
the following inequality holds true as well
\[
\lambda\ino u^2dx +
\ino|\nabla u|^2dx \geq  \frac{n^2}{4} \ino\frac{u^2}{|x|^2}dx +
\sum_{i=1}^m \int_{\xO} \frac{u^2}{|x|^2}X_1^2\ldots X_i^2 dx
 +b_m(\lambda)\ino Vu^2dx   +C \ino Wu^2dx           \ ,
\]
for $u\in C^{\infty}_{c}(\xO)$.
In case there does not exist  such a potential $W$,  we say that the potential
\[
\frac{n^2}{4} \frac{1}{|x|^2} +\sum_{i=1}^m  \frac{X_1^2\ldots X_i^2 }{|x|^2} +b_m(\lambda)V(x) \ ,
\]
is an   {\it  m--maximal potential}. Our next goal is to characterize $m$--maximal potentials. In this direction 
for $V\in\cA_m(\Omega)$  and small $r>0$ we define
\[
C_{m,r}(V) =\inf_{u\in C_{c}^{\infty}(\Omega\cap B_r)}
\frac{\lambda\int_{\Omega\cap B_r}u^2dx +
\int_{\Omega\cap B_r}|\nabla u|^2dx - \frac{n^2}{4} 
\int_{\Omega\cap B_r}\frac{u^2}{|x|^2}dx
+\frac{1}{4} \sum_{i=1}^m \int_{\xO \cap B_r } \frac{u^2}{|x|^2}X_1^2\ldots X_i^2 dx}{\int_{\Omega\cap B_r}Vu^2dx},
\]
where $X_i=X_i(|x|/3 \xk \tilde{D})$ and $\tilde{D}\geq D$.
We also define
\[
\cC^0_m(V)=\lim_{r\to 0+}C_{m,r}(V).
\]
Arguing as in Section \ref{section:impr} we can see that $\cC^0_m(V)$ is independent of the specific choice of $\lambda\geq 0$ and
$\tilde{D}\geq D$.
\begin{definition}
The potential $V\in \cA_m(\Omega)$ is $m$-subcritical if $\cC_m^0 =+\infty$.
\end{definition}
\begin{lemma}
Let $V$ be a non-negative potential satisfying
\[
\ino V^{n/2}(X_1\ldots X_{m+1})^{1-n}dx <+\infty,
\]
where $X_i=X_i(|x|/3 \xk D)$, $i=1,\ldots, m+1$. Then $V$ is an $m$-subcritical potential.
\la{lem:24}
\end{lemma}
{\bf \em Proof.} The proof is quite similar to the proof of Lemma \ref{lem:77} and makes use of Theorem \ref{g}
to establish the inequality 
\bean
&&\hspace{-3cm} \frac{\int_{\Omega\cap B_r}|\nabla u|^2dx - \frac{n^2}{4} \int_{\Omega\cap B_r}\frac{u^2}{|x|^2}dx
-\frac{1}{4}\sum_{k=1}^m \int_{\Omega\cap B_r}\frac{u^2}{|x|^2}X_1^2\ldots X_k^2 dx}
{\int_{\Omega\cap B_r}Vu^2dx} \\
&& \hspace{3cm} \geq \frac{1}{ C_n} \Big(   \int_{\Omega\cap B_r} V^{n/2}(X_1 \ldots X_{m+1})^{1-n}dx \Big)^{-\frac{2}{n}}.
\eean
The result then follows. $\hfill\Box$

As in Section \ref{section:impr} we shall also consider the following more general situation. We consider non-negative potentials
$V,W_1,W_2\in \cA_m(\Omega)$ and assume that there exist
$c>0$ and a radius $R>0$ so that
\bea
 \int_{\Omega\cap B_{R}}W_1 u^2dx +
\int_{\Omega\cap B_{R}}|\nabla u|^2dx &\geq& \frac{n^2}{4} \int_{\Omega\cap B_{R}}\frac{u^2}{|x|^2}dx  
+ \frac{1}{4} \sum_{i=1}^m \int_{\xO \cap B_r } \frac{u^2}{|x|^2}X_1^2\ldots X_i^2 dx \nonumber \\
&& + \int_{\Omega\cap B_{R}}W_2 u^2dx + c\int_{\Omega\cap B_{R}}V u^2dx 
\la{w1w2m}
\eea
for all $u\in C^{\infty}_c(\Omega\cap B_{R})$. For $0<r\leq R$ we define
\bean
& &\hspace{-6mm} C_{m,r}(W_1,W_2;V) =\\ 
&&  \hspace{-6mm} \inf_{C^{\infty}_c(\Omega\cap B_{r})}\hspace{-14pt}\frac{\int_{\Omega\cap B_{r}}\! W_1 u^2dx +\int_{\Omega\cap B_{r}}|\nabla u|^2dx - \frac{n^2}{4} \int_{\Omega\cap B_{r}}\frac{u^2}{|x|^2}dx - 
\frac{1}{4} \sum_{i=1}^m \int_{\xO \cap B_r } \frac{u^2}{|x|^2}X_1^2\ldots X_i^2 dx -
\int_{\Omega\cap B_{r}}\! W_2 u^2dx}{\int_{\Omega\cap B_{r}}V u^2dx}
\eean
and we denote
\[
\cC_m^0(W_1,W_2;V) =\lim_{r\to 0+} C_{m,r}(W_1,W_2;V).
\]
The proof of the following lemma is similar to the proof of Lemma \ref{lem:w1w2} and is omitted.
\begin{lemma}
Let $V,W_1,W_2$ be non-negative potentials in $\cA_m(\Omega)$ and assume that  there exist $R>0$ and
$c>0$ such that (\ref{w1w2m}) holds true. If in addition $W_1,W_2$ are subcritical then
$\cC_m^0(W_1,W_2;V) =\cC_m^0(V)$.
\end{lemma}

Let $\tilde{D}\geq D$ be fixed. Given $u\in C^{\infty}_c(\Omega)$ we define the function $w$ by
\bean
u(x) &=& 
 \frac{|x+\rho  e_n|^2-\rho^2}{|x|^{\frac{n}{2}}   |x+2\rho e_n|^{\frac{n}{2}}} X_1^{-\frac{1}{2}} \ldots X_m^{-\frac{1}{2}} w(x)  \\
&=:& \phi_m(x)w(x) ,  
\eean
where $X_i=X_i(|x|/3 \xk \tilde{D})$, $i=1,\ldots,m$.
Then $w\in C^{\infty}_c(\Omega)$ by our assumption that $B(-2\rho e_n,2\rho)$ is an exterior ball. 
After some computations we arrive at 
\bean
 \ino |\nabla u|^2dx &=&\ino  \phi_m^2  |\nabla w|^2dx +  n^2 \rho^2 \ino  \frac{  \phi_m^2w^2   }{|x|^{2}   |x+2\rho e_n|^{2}} dx  \\
& +& \ino \frac{\phi_m^2 w^2}{|x|^2} \bigg\{ \frac{1}{4}\sum_{k=1}^mX_1^2 \ldots X_k^2   -\frac{n}{2} \frac{ |x|^2+2\rho x_n}{|x+2\rho e_n|^2}  
\sum_{k=1}^m X_1 \ldots X_k +\frac{|x|^2}{|x|^2+2\rho x_n} \sum_{k=1}^mX_1 \ldots X_k   \bigg\}dx,
\eean
so inequality (\ref{inma}) is written
 \bea
&& \lambda\ino \phi_m^2  w^2dx + \ino  \phi_m^2 |\nabla w|^2dx +\ino \frac{\phi_m^2 w^2}{|x|^2}
\bigg\{  -\frac{n}{2} \frac{ |x|^2+2\rho x_n}{|x+2\rho e_n|^2}  
 +\frac{|x|^2}{|x|^2+2\rho x_n}  \bigg\} \Big(\sum_{k=1}^m X_1 \ldots X_k \Big) dx \nonumber \\
 & \geq & \frac{n^2}{4}\ino  \frac{ |x|^2+4\rho x_n }{|x|^2   |x+2\rho e_n|^2}   \phi_m^2 w^2dx 
+ c \ino   V \phi_m^2 w^2 dx,
\label{w1}
\eea
It is clear from the above that (\ref{w1}) is valid for all functions $w\in C^{\infty}_c(\Omega)$ and moreover,
for a fixed $\lambda \geq 0$ the best constants $c$ for inequalities (\ref{inma}) and (\ref{w1}) coincide.
That common best constant shall be denoted by $b_m=b_m(\lambda)$.

Defining
\[
Q_m(x)=    \frac{n^2}{4} \frac{ |x|^2+4\rho x_n }{|x|^2   |x+2\rho e_n|^2}    +
\frac{1}{|x|^2}
\bigg\{  \frac{n}{2} \frac{ |x|^2+2\rho x_n}{|x+2\rho e_n|^2}  
 -\frac{|x|^2}{|x|^2+2\rho x_n}  \bigg\} \Big(\sum_{k=1}^m X_1 \ldots X_k \Big) 
\]
it then follows that
\[
 C_{m,r}(V) = 
\inf_{ W^{1,2}_0( \Omega\cap B_r ; \phi_m^2)}
\hspace{-0.3cm}\frac{\lambda \int_{\Omega\cap B_r}\phi_m^2   w^2dx +
\int_{\Omega\cap B_r} \phi_m^2|\nabla w|^2dx -
 \int_{\Omega\cap B_r} Q_m  \phi_m^2 w^2dx}
{  \int_{\Omega\cap B_r} \phi_m^2 V w^2 dx}
\]
where $W^{1,2}_0( \Omega\cap B_r ; \phi_m^2)$ denotes the closure of $C^{\infty}_c(\Omega\cap B_r)$ under the norm
\[
 \int_{\Omega\cap B_r}\phi_m ^2|\nabla w|^2dx + \int_{\Omega\cap B_r} \phi_m^2 w^2dx .
\]

Similarly to Section \ref{section:impr} we shall use a simpler way for expressing $\cC^0_m(V)=\lim_{r\to 0+}C_{m,r}(V)$. For this we define
\[
C_{m,r}(V;\phi_m^2) =\inf_{W^{1,2}_0( \Omega\cap B_r ; \phi_m^2)}
\frac{ \int_{\Omega\cap B_r}  \phi_m^2|\nabla w|^2dx}{ \int_{\Omega\cap B_r}  V \phi_m^2  w^2dx}.
\]
and
\[
\cC^0_m(V;\phi_m^2) =\lim_{r\to 0+} C_{m,r}(V;\phi_m^2).
\]

\begin{lemma}
Let $V$ be a non-negative potential in $\cA_m(\Omega)$. Then
$\cC^0_m(V) =\cC^0_m(V;\phi_m^2)$.
\end{lemma}
{\bf \em Proof.} The proof is quite similar to the proof of Lemma \ref{lem:fr_ir}. In particular we make use of the fact that
\[
\int_{\Omega} \Big(\frac{X_1}{|x|^2+2x_n}\Big)^{n/2} (X_1\ldots X_{m+1})^{1-n} < +\infty ,
\]
from which it easily follows that $|Q_m|$ is an $m$-subcritical potential. We omit further details. $\hfill\Box$

One important consequence of $m$-subcriticality is the following compactness property whose proof is similar to 
that of Lemma \ref{compactness} and is therefore omitted.
\begin{lemma}
Assume that the positive potential $V\in\cA_m(\Omega)$ is $m$-subcritical. Then for any sequence $(w_k)$ which is bounded in $W^{1,2}_0(\Omega ;\phi_m^2)$ there exists a subsequence, also denoted by $(w_k)$, and a function $w_0\in W^{1,2}_0(\Omega ;\phi_m^2)$ so that
\bean
\ia && w_k \rightharpoonup w_0 \mbox{ in } W^{1,2}_0(\Omega ;\phi_m^2)\\
\ib && \int_{\Omega} V \phi_m^2(w_k-w_0)^2dx \longrightarrow 0 .
\eean
\end{lemma}

We can now state and prove the main theorems of this section.
\begin{theorem}
Let $V\in\cA_m(\Omega)$ and $\lambda\geq 0$ be given and let $b_m(\lambda)$ be the best constant for the inequality
\be
 \lambda\ino \phi_m^2  w^2dx + \ino \phi_m^2 |\nabla w|^2dx \geq  \ino  Q_m  \phi_m^2 w^2dx  + c \ino   V \phi_m^2 w^2 dx,
 \qquad w\in  W^{1,2}_0(\Omega;\phi_m^2),
\la{eq:insta}
\ee
where
\[
 \phi_m(x)=
\frac{|x+\rho  e_n|^2-\rho^2}{|x|^{\frac{n}{2}}   |x+2\rho e_n|^{\frac{n}{2}}} X_1^{-\frac{1}{2}} \ldots X_m^{-\frac{1}{2}}
\]
and
\[
Q_m(x)=    \frac{n^2}{4} \frac{ |x|^2+4\rho x_n }{|x|^2   |x+2\rho e_n|^2}    +
\frac{1}{|x|^2}
\bigg\{  \frac{n}{2} \frac{ |x|^2+2\rho x_n}{|x+2\rho e_n|^2}  
 -\frac{|x|^2}{|x|^2+2\rho x_n}  \bigg\} \Big(\sum_{k=1}^m X_1 \ldots X_k \Big) ,
\]
with $X_i=X_i(|x|/3  \xk \tilde{D})$, $\tilde{D}\geq D$.
If in addition
\[
b_m(\lambda)<\cC_m^0
\]
then the best constant $b_m(\lambda)$ in (\ref{eq:insta}) is realized by a function $w_0\in W^{1,2}_0(\Omega;\phi_m^2)$. In particular the best constant $b_m(\lambda)$ is realized if $V$ is an $m$-subcritical potential.
\label{gera:thm:exmin1}
\end{theorem}
{\bf \em Proof.} The proof is similar to the proof of Theorem \ref{gera:thm:exmin} so we shall only give a sketch of the proof. What is important for our argument is that the potential $Q_m$ is a subcritical potential.
We consider a minimizing sequence $(w_k)$ for (\ref{eq:insta}) and without loss of generality we assume that
\[
\ino  Q_m \phi_m^2  w_k^2dx  +  b_m(\lambda) \ino V \phi_m^2  w_k^2 dx
=1, \qquad   \lambda\ino \phi_m^2  w_k^2dx + \ino  \phi_m^2 |\nabla w_k|^2dx  \longrightarrow 1, \; \mbox{ as }k\to \infty.
\]
Since $(w_k)$ is bounded in $W^{1,2}_0(\Omega ; \phi_m^2)$, it has a subsequence, which we assume is $(w_k)$ itself, which converges weakly to some $w_0\in W^{1,2}_0(\Omega ; \phi_m^2)$. We define $v_k=w_k-w_0$.

We consider a small enough $r>0$ so that $C_{m,r}(V;\phi_m^2)>b_m(\lambda)$ and a smooth cut-off function $\psi$
such that $\psi=1$ in $B_{r/2}$ and $\psi=0$ outside $B_r$. Arguing as in the proof of Theorem \ref{gera:thm:exmin} we obtain
\be
\int_{\Omega}\phi_m^2 |\nabla v_k|^2 dx \geq  C_{m,r}(V;\phi_m^2) \int_{\Omega}V \phi_m^2  v_k^2 dx +o(1) 
\la{equm}
\ee
and also 
\be
\ino Q_m\phi_m^2  w_0^2dx  +  b_m(\lambda) \ino V \phi_m^2  w_0^2 dx + b_m(\lambda) \ino V \phi_m^2  v_k^2 dx  =1 +o(1)
\la{norm1a}
\ee
and
\be
\lambda \ino \phi_m^2  w_0^2dx + \ino  \phi_m^2 |\nabla w_0|^2dx + \ino  \phi_m^2 |\nabla v_k|^2dx = 1+o(1).
\la{norm2a}
\ee
From (\ref{equm}) and (\ref{norm1a}) we obtain
\[
\int_{\Omega}\phi_m^2 |\nabla v_k|^2 dx \geq \frac{ C_{m,r}(V;\phi_m^2)}{b_m(\lambda)}\Big(  1-    \ino Q_m \phi_m^2  w_0^2dx  +  b_m(\lambda) \ino V \phi_m^2  w_0^2 dx        \Big) +o(1).
\]
Writing (\ref{eq:insta}) for $w=w_0$ we obtain from (\ref{norm2a}) that
\[
\ino  \phi_m^2 |\nabla v_k|^2dx \leq 1 -\int_{\Omega}Q_m \phi_m^2w_0^2dx -b_m(\lambda)\ino V\phi_m^2w_0^2dx +o(1) 
\]
and arguing as before we conclude that
\[
 \ino Q_m \phi_m^2  w_0^2dx  +  b_m(\lambda) \ino V \phi_m^2 V w_0^2 dx =1,
\]
that is $w_0$ is a minimizer for (\ref{eq:insta}). 
$\hfill\Box$

Finally, the next is a direct consequence of Theorem \ref{gera:thm:exmin1}.
\begin{theorem}
\label{h}
Let $n \geq 3$.  $V$ is a non-negative potential in $\cA_m(\Omega)$ and $\tilde{D}\geq D$.
(a) Let $\lambda\geq 0$ and $b_m(\lambda)>0$ be the best constant for the following inequality
\be
\lambda\ino u^2dx +
\ino|\nabla u|^2dx \geq  \frac{n^2}{4} \ino\frac{u^2}{|x|^2}dx + \frac{1}{4} \sum_{k=1}^m \int_{\xO} \frac{u^2}{|x|^2}X_1^2\ldots X_k^2 dx
+ b_m(\lambda)\ino Vu^2dx, \quad u\in C^{\infty}_{c}(\Omega),
\la{in60}
\ee
where $X_i=X_i(|x|/3 \xk \tilde{D})$. If in addition
\[
b_m(\lambda)<\cC_m^0 \ ,
\]
then the potential $n^2/4|x|^2 + 1/4|x|^2 \sum_{k=1}^m X_1^2\ldots X_k^2  +  b_m(\lambda)V(x)$ is a maximal potential, that is, inequality
(\ref{in60}) cannot be improved by adding a non-negative potential $W$ in the RHS. \nl
(b) If $V$ is an $m$-subcritical potential then there exist $\lambda\geq 0$ and a best constant $b_m(\lambda)>0$ such that (\ref{in60}) is true. Moreover the potential $n^2/4|x|^2 +  1/4|x|^2 \sum_{k=1}^m X_1^2\ldots X_k^2 + b_m(\lambda)V(x)$ is a maximal potential.
\end{theorem}


\section{Maximal potentials in finite cones}

In the previous two sections we characterized maximal potentials in bounded domains satisfying the exterior and interior ball condition. Analogous results  also hold true in the case of finite cones.

Let $\ccC_1$  be the cone  determined by the domain $\Sigma\subset S^{n-1}$ as defined in Section \ref{section:cone}; more
generally  we set $\ccC_r:=\ccC \cap B_r$.
In this subsection we shall initially  be interested in the question of characterizing maximal potentials
for improved versions of inequality (\ref{cone}).
\begin{definition}
A non-negative potential $V\in L^{n/2}_{\loc}(\overline{\ccC_1}\setminus\{0\})$ is called {\em admissible} if there exists
$c>0$ such that
\be
\int_{\ccC_1}|\nabla u|^2dx \geq  \left(  \Big(\frac{n-2}{2}\Big)^2  + \mu_1(\xS) \right) \int_{\ccC_1} \frac{u^2}{|x|^2} dx
+c  \int_{\ccC_1} Vu^2dx, \quad u\in C^{\infty}_{c}(\ccC_1).
\la{1234c}
\ee
We denote by $\cA(\ccC_1)$ the class of all admissible potentials. 
\end{definition}
 
Once again there is a big variety of admissible potentials. For example if $V$ satisfies 
$\int_{\ccC_1} V^{n/2}X_1^{1-n}dx <+\infty$
where $X_1=X_1(|x|)$, then $V$ is admissible by Theorem \ref{i}.

Given $V\in\cA(\ccC_1)$ and $r\in (0,1)$ we define
\[
C_r(V) =\inf_{u\in C^{\infty}_c(\ccC_r)}
\frac{\int_{\ccC_r}|\nabla u|^2dx - \left(  \Big(\frac{n-2}{2}\Big)^2  + \mu_1(\xS) \right)   
\int_{\ccC_r}\frac{u^2}{|x|^2}dx}{\int_{\ccC_r}Vu^2dx}
\]
and
\[
\cC^0(V)=\lim_{r\to 0+}C_r(V).
\]

\begin{definition}
We say that the potential $V\in \cA(\ccC_1)$ is subcritical if $\cC^0(V) =+\infty$.
\end{definition}

The analogue of Lemma \ref{lem:77} is the following
\begin{lemma}
Let $V$ be a non-negative potential satisfying
\[
\int_{\ccC_1} V^{n/2}X_1^{1-n}dx <+\infty,
\]
where $X_1=X_1(|x|)$. Then $V$ is a subcritical potential.
\end{lemma}

Given $u\in C^{\infty}_c(\ccC_1)$ we define the function $w$ by
\be
u(x) = |x|^{-\frac{n-2}{2}}\phi_1(\omega) w(x) =:\psi(x) w(x),
\la{uwc}
\ee
where $\omega:=\frac{x}{|x|}  \in \Sigma$  and  $\phi_1(\omega)$, is the first eigenfunction of the Dirichlet Laplacian in $\Sigma$. 
After some computations we obtain
\[
\int_{\ccC_1}|\nabla u|^2dx =\int_{\ccC_1}  \psi^2  |\nabla w|^2dx + \left(  \Big(\frac{n-2}{2}\Big)^2  + \mu_1(\xS) \right)   \int_{\ccC_1}  \psi^2 w^2dx .
\]
As usual we denote by 
$W^{1,2}_0(\ccC_1  ; \psi^2)$ the closure of $C^{\infty}_c(\ccC_1)$ under the norm
\[
\int_{\ccC_1}\psi^2|\nabla w|^2dx +
 \int_{\ccC_1} \psi^2  w^2dx .
\]
It is easily seen that inequality (\ref{1234c}) 
 under the change of variables (\ref{uwc}) is equivalent to
\[
\int_{\ccC_1} \psi ^2|\nabla w|^2dx \geq 
  c \int_{\ccC_1} V \psi^2 w^2dx , \qquad w\in W^{1,2}_0(\ccC_1;\psi^2).
\]
The analogues of Theorems \ref{gera:thm:exmin} and \ref{d} read as follows: 

\begin{theorem}
Let $V\in\cA(\ccC_1)$  and suppose that $b$ is the best constant for the inequality
\be
\int_{\ccC_1} \psi ^2|\nabla w|^2dx \geq 
  b \int_{\ccC_1} V \psi^2 w^2dx , \qquad w\in W^{1,2}_0(\ccC_1;\psi^2).
\la{eq:w:no1c}
\ee
If in addition $b<\cC^0(V)$ 
then the best constant $b$ in (\ref{eq:w:no1c}) is realized by a function $w_0\in W^{1,2}_0(\ccC_1;\psi^2)$. In particular the best constant $b$ is realized if the potential $V$ is subcritical.
\end{theorem}

\begin{theorem}
Let $V$ be a non-negative potential in $\cA(\ccC_1)$. \nl
 (a) Let $b>0$ be the best constant in the following inequality
\be
\int_{\ccC_1}|\nabla u|^2dx \geq  \bigg(  \Big(\frac{n-2}{2}\Big)^2  + \mu_1(\xS) \bigg) 
\int_{\ccC_1}\frac{u^2}{|x|^2}dx +b \int_{\ccC_1} Vu^2dx,
 \quad u\in  C^{\infty}_{c}(\ccC_1) \ .
\la{eq:thm:123:c}
\ee
 If in addition $b<\cC^0(V)$
then the potential
$\big[  \big(\frac{n-2}{2}\big)^2  + \mu_1(\xS) \big]   |x|^{-2} +b V(x)$
is a maximal potential, that is inequality
(\ref{eq:thm:123:c}) cannot be improved by adding a non-negative potential $W$ in the RHS. \nl
\noindent (b) If $V$ is a subcritical potential then there exists a best constant $b>0$ such that (\ref{eq:thm:123:c}) is true. Moreover the potential
$\big[  \big(\frac{n-2}{2}\big)^2  + \mu_1(\xS) \big]   |x|^{-2}   +bV(x)$ is a maximal potential.
\end{theorem}

The proofs of these theorems are quite similar and slightly simpler to the proofs of of Theorems \ref{gera:thm:exmin} and \ref{d}.

In analogy to the results of Section \ref{sub:log:impr}  we have  similar theorems for the improved Hardy inequality involving logarithmic corrections. In particular we have

\begin{definition}
A non-negative potential $V\in L^{n/2}_{\loc}(\overline{\ccC_1}\setminus\{0\})$ is called {\em m--admissible} if there exists
$c>0$ such that for $ u\in C^{\infty}_{c}(\ccC_1)$, there holds
\be
\int_{\ccC_1}|\nabla u|^2dx \geq  \left(  \Big(\frac{n-2}{2}\Big)^2  + \mu_1(\xS) \right) \int_{\ccC_1} \frac{u^2}{|x|^2} dx
+\frac{1}{4}\sum_{i=1}^m
\int_{\ccC_{1}}\frac{u^2}{|x|^2}X_1^2\ldots X_i^2 dx  
+c  \int_{\ccC_1} Vu^2dx \ .
\la{1234cc}
\ee
\end{definition}
We denote by $\cA_m(\ccC_1)$ the class of all $m$--admissible potentials.

Given $V\in\cA_m(\ccC_1)$ and $r\in (0,1)$ we define
\[
C_{m,r}(V) =\inf_{u\in C^{\infty}_c(\ccC_r)}
\frac{\int_{\ccC_r}|\nabla u|^2dx - \bigg(  \Big(\frac{n-2}{2}\Big)^2  + \mu_1(\xS) \bigg)   
\int_{\ccC_r}\frac{u^2}{|x|^2}dx- \frac{1}{4}\sum_{i=1}^m
\int_{\ccC_{r}}\frac{u^2}{|x|^2}X_1^2\ldots X_i^2 dx  }{\int_{\ccC_r}Vu^2dx}
\]
and
\[
\cC^0_m(V)=\lim_{r\to 0+}C_{m,r}(V).
\]

\begin{definition}
We say that the potential $V\in \cA_m(\ccC_1)$ is $m$-subcritical if $\cC^0_m(V) =+\infty$.
\end{definition}

Changing variables by 
\[
u(x) =|x|^{-\frac{n-2}{2}}\phi_1(\omega)X_1^{-\frac{1}{2}}\ldots  X_m^{-\frac{1}{2}} w(x)=:\psi_m(x) w(x) \ ,
\]
inequality (\ref{1234cc}) is equivalent to
\[
\int_{\ccC_1} \psi_m^2|\nabla w|^2dx \geq 
  c \int_{\ccC_1} V \psi^2_m w^2dx , \qquad w\in W^{1,2}_0(\ccC_1;\psi^2_m).
\]

Now the analogues of Theorems \ref{gera:thm:exmin1} and \ref{h} are as follows

\begin{theorem}
Let $V\in\cA_m(\ccC_1)$ and let $b_m$ be the best constant for the inequality
\be
  \int_{\ccC_1}  \psi_m^2 |\nabla w|^2dx \geq  b_m \int_{\ccC_1}  V \psi_m^2 w^2 dx,
 \qquad w\in  W^{1,2}_0(\Omega;\psi_m^2),
\la{eq:instac}
\ee
If in addition $b_m<\cC^0_m(V)$ 
then the best constant $b_m$ in (\ref{eq:instac})  is realized   in  $W^{1,2}_0(\ccC_1;\psi^2_m)$. In particular the best constant $b_m$ is realized if the potential $V$ is $m$--subcritical.
\end{theorem}

\begin{theorem}
Let $V$ be a non-negative potential in $\cA_m(\ccC_1)$. \nl
 (a) Let $b_m>0$ be the best constant in the following inequality
\be
\int_{\ccC_1}|\nabla u|^2dx \geq  \bigg(  \Big(\frac{n-2}{2}\Big)^2  + \mu_1(\xS) \bigg)\int_{\ccC_1}\frac{u^2}{|x|^2}dx +\frac{1}{4}\sum_{i=1}^m   \int_{\ccC_{1}}\frac{u^2}{|x|^2}X_1^2\ldots X_i^2 dx +
+b_m \int_{\ccC_1} Vu^2dx, 
\la{eq:thm:123:cc}
\ee
where $u\in  C^{\infty}_{c}(\ccC_1)$.  \nl
 If in addition $b_m<\cC^0_m(V)$
then the potential
$\big[  \big(\frac{n-2}{2}\big)^2  + \mu_1(\xS)    
 + \frac{1}{4}\sum_{i=1}^m
X_1^2\ldots X_i^2  \big] |x|^{-2} +b_m V(x)$
is a maximal potential, that is inequality
(\ref{eq:thm:123:cc}) cannot be improved by adding a non-negative potential $W$ in the RHS. \nl
\noindent (b) If $V$ is a subcritical potential then there exists a best constant $b_m>0$ such that (\ref{eq:thm:123:cc}) is true. Moreover the potential
$\big[  \big(\frac{n-2}{2}\big)^2  + \mu_1(\xS) + \frac{1}{4}\sum_{i=1}^m
X_1^2\ldots X_i^2   \big]   |x|^{-2}   +b_mV(x)$ is a maximal potential.
\end{theorem}


{\bf Remark.} All the above results involve a single point singularity on the boundary. Similar results however can be obtained when there
are multiple singularities. For instance we have the following result

\begin{theorem}
\label{th72}
Assume that $\Omega \subset \R^n$, $n \geq 3$, is a bounded domain that satisfies an exterior ball condition at each of the points $a_1,\ldots,a_m\in \partial\Omega$.
Then there exist a positive constant $c=c(n,m)$ depending only on $n$ and $m$ and a positive constant $\lambda$ such that
\[
\lambda \int_{\Omega}u^2 dx + \int_{\Omega}|  \nabla u|^2 dx \geq \frac{n^2}{4}\sum_{k=1}^m \ino\frac{u^2}{|x-a_k|^2}dx
+c \Big( \ino |u|^{\frac{2n}{n-2}}W dx\Big)^{\frac{n-2}{n}} \ , ~~~~~u\in C^{\infty}_c(\Omega) \ ,
\]
where 
\[ 
 W=W(x):=\prod_{k=1}^{m} X_1^{\frac{2n-2}{n-2}}(\frac{|x-a_k|}{D}),~~~~~{\rm and }~~~~ D:=  \max_{k=1,\ldots, m}~ \sup_{x \in \xO} |x-a_k|  \ .
 \]
\end{theorem}
The proof uses ideas that we have used so far in connection with standard partition of unity arguments; we omit further details.



\medskip

{ \bf Acknowledgement} We would like to thank the referees for their comments which led to a considerable improvement of the presentation
of this work. AT acknowledges partial support by ELKE grant, University of Crete.



\end{document}